\documentclass[12pt,reqno,twoside]{amsart}
\usepackage{amsfonts,amsmath,amssymb}
\usepackage{mathrsfs,mathtools}
\usepackage{enumerate}
\usepackage{hyperref}
\usepackage{esint}
\usepackage{graphicx,subfigure}
\usepackage{bm}
\usepackage{commath}
\usepackage{esint}
\usepackage{placeins}
\usepackage{flafter}
\usepackage{tikz}
\usepackage[textsize=small]{todonotes}
\usepackage[dvips,bottom=1.4in,right=1in,top=1in, left=1in]{geometry}
\usepackage[latin1]{inputenc}

\usepackage{caption}

\DeclareMathAlphabet{\mathpzc}{OT1}{pzc}{m}{it}
\numberwithin{equation}{section}

\makeatletter
\def\eqnarray{\stepcounter{equation}\let\@currentlabel=\theequation
\global\@eqnswtrue
\tabskip\@centering\let\\=\@eqncr
$$\halign to \displaywidth\bgroup\hfil\global\@eqcnt\z@
  $\displaystyle\tabskip\z@{##}$&\global\@eqcnt\@ne
  \hfil$\displaystyle{{}##{}}$\hfil
  &\global\@eqcnt\tw@ $\displaystyle{##}$\hfil
  \tabskip\@centering&\llap{##}\tabskip\z@\cr}
\def\endeqnarray{\@@eqncr\egroup
      \global\advance\c@equation\m@ne$$\global\@ignoretrue}

\numberwithin{equation}{section}

\def\pmb#1{\setbox0=\hbox{$#1$}%
             \kern-.027em\copy0\kern-\wd0
             \kern+.009em\copy0\kern-\wd0
             \kern+.009em\copy0\kern-\wd0
             \kern+.009em\copy0\kern-\wd0
             \kern+.009em\copy0\kern-\wd0
             \kern+.009em\copy0\kern-\wd0
             \kern+.009em\copy0\kern-\wd0
             \kern-.045em\raise+.012em\copy0\kern-\wd0
             \kern+.009em\raise+.012em\copy0\kern-\wd0
             \kern+.009em\raise+.012em\copy0\kern-\wd0
             \kern+.009em\raise-.012em\copy0\kern-\wd0
             \kern+.009em\raise-.012em\copy0\kern-\wd0
             \kern-.018em\copy0\kern-\wd0\raise-.012em\box0}
\def\Pmb#1{\setbox0=\hbox{$#1$}%
             \kern-.033em\copy0\kern-\wd0
             \kern+.011em\copy0\kern-\wd0
             \kern+.011em\copy0\kern-\wd0
             \kern+.011em\copy0\kern-\wd0
             \kern+.011em\copy0\kern-\wd0
             \kern+.011em\copy0\kern-\wd0
             \kern+.011em\copy0\kern-\wd0
             \kern-.055em\raise+.015em\copy0\kern-\wd0
             \kern+.011em\raise+.015em\copy0\kern-\wd0
             \kern+.011em\raise+.015em\copy0\kern-\wd0
             \kern+.011em\raise-.015em\copy0\kern-\wd0
             \kern+.011em\raise-.015em\copy0\kern-\wd0
             \kern-.022em\copy0\kern-\wd0\raise-.015em\box0}
\def\ms{\medskip}

\def\ub{\underbar}
\def\noi{\noindent}

\def\ptilde{\tilde{p}}
\def\fvec{{\bf f}}
\def\gvec{{\bf g}}

\def\rvec{{\bf r}}
\def\svec{{\bf s}}

\def\uvec{{\bf u}}
\def\vvec{{\bf v}}

\def\xvec{{\bf x}}

\def\dt{{\Delta t}}

\def\Grd{\nabla}
\def\Div{\nabla \cdot}

\vfill

\title[DETERMINISTIC PATHOGEN TRANSMISSION MODEL]{A DETERMINISTIC PATHOGEN TRANSMISSION MODEL BASED ON 
HIGH-FIDELITY PHYSICS}
\date{\today}

\thanks{
H. Antil is partially supported by NSF grants DMS-2110263, DMS-1913004,
the Air Force Office of Scientific Research (AFOSR) under
Award NO: FA9550-19-1-0036, and Department of Navy, Naval PostGraduate
School under Award NO: N00244-20-1-0005. \\
E. O\~nate and S. Idelsohn acknowledge financial support from the 
project PARAFLUIDS (PID2019-104528RB-I00) of the National Research 
Plan of the Spanish Government, from the CERCA programme of 
the Generalitat de Catalunya, and from the Spanish Ministry of 
Economy and Competitiveness, through the `Severo Ochoa Programme 
for Centres of Excellence in R\&D' (CEX2018-000797-S).
}

\author[R.~L\"ohner]{Rainald L\"ohner}
\address{R. L\"ohner. Center for Computational Fluid Dynamics, College of Science, George Mason University, Fairfax, VA 22030-4444, USA.}
\email{rlohner@gmu.edu}
\author[H. Antil]{Harbir Antil}
\address{H. Antil. The Center for Mathematics and Artificial Intelligence
(CMAI) and Department of Mathematical Sciences, George Mason University,
Fairfax, VA 22030, USA.}
\email{hantil@gmu.edu}
\author[J.M~Gimenez]{Juan Marcelo Gimenez}
\address{J.M. Gimenez. CIMNE, International Center for Numerical Methods in
            Engineering, Barcelona, Spain.}
\email{jmgimenez@cimne.upc.edu}
\author[S.~Idelsohn]{Sergio Idelsohn}
\address{S. Idelsohn.  ICREA, Catalan Institution for Research and Advanced Studies
            and CIMNE, International Center for Numerical Methods in
            Engineering, Barcelona, Spain.}
\email{sergio@cimne.upc.edu}
\author[E.~O\~nate]{Eugenio O\~nate}
\address{E. O\~nate. CIMNE, International Center for Numerical Methods in
            Engineering, Barcelona, Spain.}
\email{onate@cimne.upc.edu}

\begin{document}

\begin{abstract}	
A deterministic pathogen transmission model based on high-fidelity physics
has been developed. The model combines computational fluid dynamics and
computational crowd dynamics in order to be able to provide accurate
tracing of viral matter that is exhaled, transmitted and inhaled via 
aerosols. The examples
shown indicate that even with modest computing resources, the propagation
and transmission of viral matter can be simulated for relatively large 
areas with thousands of square meters, hundreds of pedestrians and several
minutes of physical time. The results obtained and insights gained from 
these simulations can be used to inform global pandemic propagation models, 
increasing substantially their accuracy.
\end{abstract}

\keywords{Pathogen Transmission, Viral Transmission,
      Pathogen Mitigation,
      Finite Elements, 
      Computational Fluid Dynamics, 
      Computational Crowd Dynamics}

\maketitle
\tableofcontents
	
\section{Introduction}
\label{intro}
\begin{figure}
	\includegraphics[width=10.0cm]{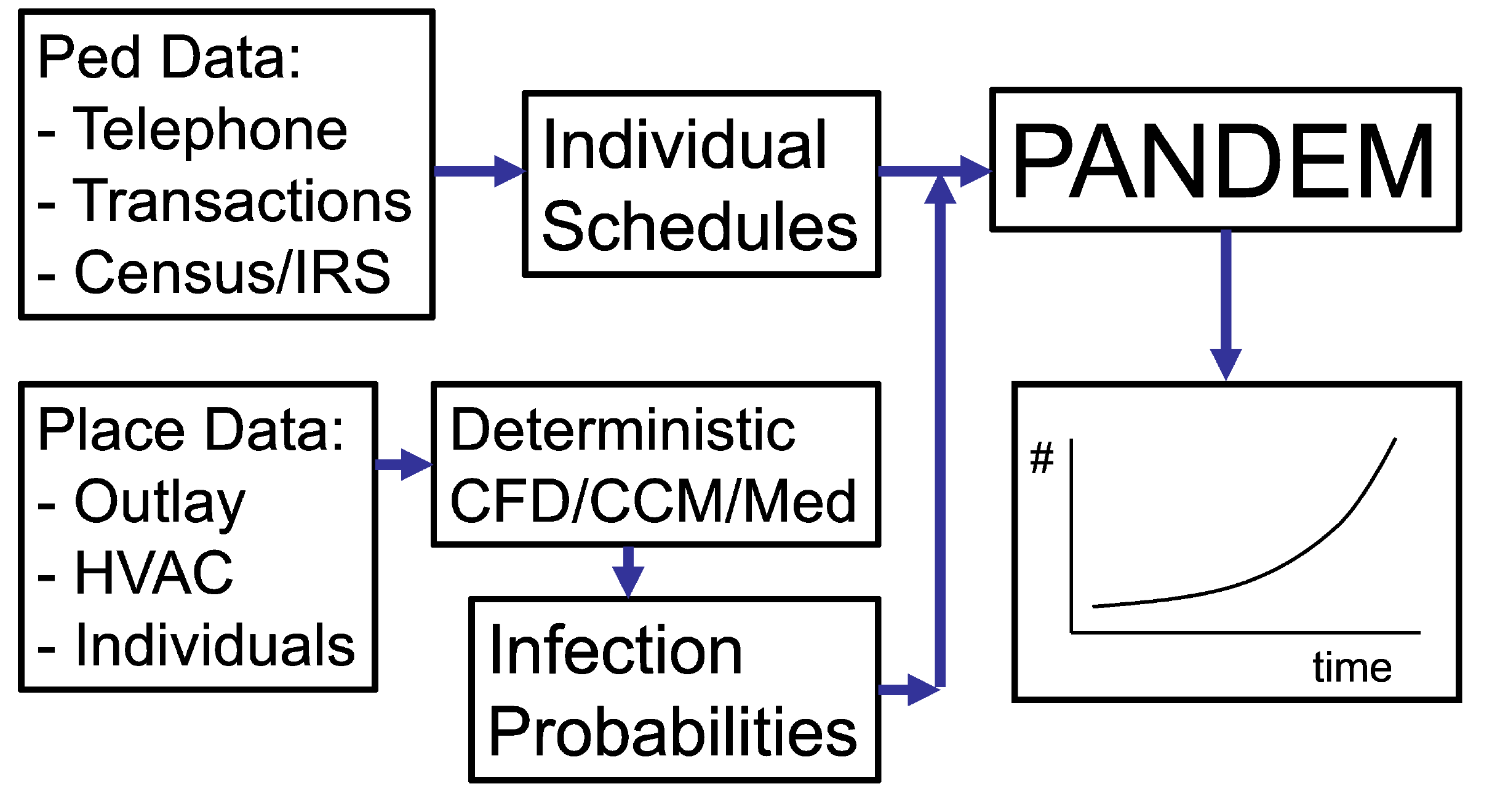}
	\caption{Global Pandemic Simulation Codes}
	\label{f:pendem}
\end{figure}

Many public health decisions are being made based on results from global
pandemic simulation models as shown schematically in Figure~\ref{f:pendem}.
The whereabouts of people -~and hence their proximity~- during a
typical day can be obtained from a variety of sources such 
as telephone locator records and transactional records. The accuracy
of this data is relatively high, allowing to place individuals
within a few meters. Advances in both compute power as well as
algorithms and optimal data structures have enabled to track and
compute the interaction of billions of individuals during several
days on desktop workstations in a matter of hours. Thus, statistical
runs are easy to perform on larger clusters.  \\

\par \noi
The highest uncertainty of these types of models is the transmission
rate in given settings. It is easy to see that climate 
(heating, cooling, ventilation), customs (hugging, kissing, proximity)
and population behaviour (masking) can lead to vastly different 
transmission rates for the same setting (public office, classroom, 
meeting room, factory, airport, train station, cinema, theater,
stadium, restaurant, etc.). Particularly for those pathogens that are
very transmissible via the aerosol route (e.g. measles, Covid-19)
the variability can be so large that the accuracy of these global
pandemic simulation models can be called into question. \\
It is here that simulations based on high-fidelity physics of flows 
with aerosols together with pedestrian/crowd dynamics can be used
to achieve a higher accuracy in the transmission rates, thus improving 
the accuracy of global simulation models. \\
This would no doubt entail a formidable undertaking. In principle,
every venue would have to be analyzed. And this implies knowing the
outlay and dimensions of the place, the heating, ventilation and
air conditioning (HVAC) information, the 
possible external climate factors, the flow of people, etc.
Given the variability of any of these factors, each one of these 
venues would in turn require a number of runs. However, the
insight gained could be used to a)~improve the airflow and ventilation
so as to minimize pathogen transmission, b)~modify the movement of
people so as to minimize pathogen transmission, and c)~increase 
the awareness of people frequenting these venues via visualization 
of pathogen movement. \\ 
On the other hand, the variability in the number of places frequented
is also not infinite. Many crowded places belong to chains or 
franchises, and are therefore 'standardized' in their outlay.
Typical venues that fall under this category are fast-food restaurants.
The flow here is predominantly from the customer/public area towards the
kitchen/cooking area in order avoid odours. This implies that the
placement (e.g. queuing) and movement of people are key factors that
will affect transmission rates. Similar places that also have a 
'standardized outlay' are office and government buildings. \\
The present paper presents first steps in this direction. As stated
before, an obvious necessity is the development of a reliable,
accurate and scalable deterministic high-fidelity pathogen 
transmission model that incorporates the proper physics of flows,
pedestrian motion and the ex/inhalation of pathogens.

\section{Modeling Pathogen Transmission and Infection}
\label{modelingpathogentransmission}
Pathogen, or more generally, disease transmission models have a long
history in both health sciences and applied mathematical modeling
since the seminal work of Wells \cite{Wel55}. The probability
of a person getting infected through the airborne route depends on
\begin{itemize}
\item[-] That person's exposure to the pathogen/virus (the dose) 
\cite{Fra97,Bro17} and
\item[-] The probability of getting infected given that level of
exposure \cite{Sze10}.
\end{itemize}
The first term depends on a variety of factors such as pathogens shed by
infective persons, exposure time, air circulation patterns,
etc. \cite{Lou67,Tan06,Xie07,Lin10,Teu10,Sze10,Gup11b,Joh11,Tan11,Gup12,Gup12,Hal12,Lin12,Tan12,Mil13,Tan13,Wei16,Asa19,Asa20a,Asa20b,WHO20,To20}.
In aggregate models, which include the vast
majority of models, the specific locations of susceptible persons are
not explicitly considered. Rather, one makes assumptions so that these
locations do not have to be specified. For example, the viral particles
may be assumed to be uniformly distributed, in which case the specific
locations of persons are not required.  More sophisticated models
can account for spatial heterogeneity by dividing the space into
multiple zones with different distributions of viral particles in each
zone \cite{Arm07}. This requires knowing only the exposure times
of persons in each zone, and not their specific positions. In
individualized models, which are relatively rare, the specific
positions of persons are required. For example, Gupta et al. \cite{Gup12}
and L\"ohner et al. \cite{Loh21a} consider specific seats in which 
passengers in a plane are seated and examine their exposure using 
Computational Fluid Dynamics (CFD) simulations. \\
The second term - the relationship between exposure and infection
probability - can be deterministic or probabilistic \cite{Sze10,Gup12},
with the latter being much more common. In a deterministic model, a
person is considered infected if that person inhales more viral
particles than a limit, called that person's tolerance dose. In a
probabilistic model, on the other hand, the infection probability
depends on the extent of exposure. Models typically use variants of
one of the following two approaches. In both approaches, infection
probability is given by

$$ p = 1 - e^{- c \cdot e_m}  ~~, \eqno(2.1) $$

\noi
where $c$ is a constant and $e_m$ the so-called exposure metric.
In the Wells-Riley approach \cite{Wel55}, the exposure is expressed 
in terms of
an abstract `quantum' of infection, whose relative values can be
computed for different scenarios, and model parameters determined
by fitting against empirical data. In a dose-response model, on the
other hand, the exposure metric reflects the actual number of viral
particles inhaled by the susceptible person, with detailed mechanisms
for computing this value.

\subsection{Relationship to Previous Work}
Conventional models are aggregate, and do not consider the specific
locations of individuals. Consequently, they cannot account for
fine-scaled spatial heterogeneity. Individual models, such
as \cite{Gup12,Loh21}, do consider positions of individuals. 
However, they
mostly consider situations where the positions of persons are fixed. This
is inadequate for understanding risk associated when people move in
a crowd. Namilae et al. \cite{Nam17} consider movement of people in
a plane using pedestrian dynamics. However, that work does not
account for movement of viral particles through the air. Instead,
it considers infection risk based on contacts between persons in
the crowd. L\"ohner al al. \cite{Loh21} consider movement of people 
in a fully coupled setting. \\

\section{Requirements for Modeling Pathogen Propagation, Transmission 
and Mitigation}
\label{requirements}
Taking into account all the information stated before, one can see
that in order to arrive at advanced numerical models to compute with
high fidelity pathogen propagation, transmission and mitigation, the
following capabilities are required:
\begin{itemize}
\item[-] Physical modeling of sneezing/coughing (exit velocities and
temperature, number and distribution of particles, ...);
\item[-] Physical modeling of aerosol propagation (flows with 
particles in an environment with moving pedestrians, geometric 
fidelity of the built environment, HVAC boundary conditions, ...);
\item[-] Modeling of pedestrian motion (movement, proximity, ...);
\item[-] Monitoring of pathogens exhaled and inhaled.
\end{itemize}
\noindent
These in turn will enable the generation of the four essential 
pieces of information required:
\begin{itemize}
\item[-] The generation of pathogen (e.g. viral) loads;
\item[-] The movement (advection, diffusion) of pathogen loads;
\item[-] The location and movement of pedestrians exhaling pathogens;
\item[-] The location and movement of pedestrians inhaling pathogens.
\end{itemize}
In the sequel, we will consider each one of these in turn.
One should state from the outset that all of these quantities can vary
greatly, so that any kind of model will have to be run repeatedly
in order to obtain proper statistics.

\section{Physical Modeling of Aerosol Propagation}
\label{physmodelaerosols}
The question often arises whether for pathogen transmission the
motion of individual droplets needs to be computed. The larger
droplets of diameters $O(1.0~mm)$ tend to fall ballistically.
The smaller ones, with diameters $d<O(0.1~mm)$, tend to slow 
down immediately and adjust to the velocity of the surrounding air. 
Furthermore, they also evaporate quickly.
If one considers the motion of a
water particle with an initial velocity of $v_i=1~m/sec$ into
quiescent air, and the usual
values of $\rho_p=1000~kg/m^3, \rho_{air}=1.2~kg/m^3,
\mu_{air}=1.85\cdot~10^{-5}~kg/m/s$, one can obtain the
distance and time
to rest, where `rest' in this case is assumed as $v_r=0.01~v_i$.
These values have been tabulated in Table~\ref{tab:Distimerest}. One can see that for
diameters below $O(0.1~mm)$ the time and distance required for
a particle to adjust to the velocity of the surrounding air is so
low that for these aerosol particles one can neglect the
air-particle interaction. Therefore, one can treat these
aerosol particles via a transport equation that advects and diffuses
the particle concentration in space and time.

\begin{table}[htbp]
\begin{center}
\caption{Distance and Time to Rest for Water Particles in Air
With An Initial Velocity of 1m/sec}
\label{tab:Distimerest}
\begin{tabular}{c|c|c}
\hline
Diameter [mm] & distance to rest [m] & time to rest [sec] \\
\hline
1.00E-01      & 2.27E-02      & 1.20E-01 \\
1.00E-02      & 2.79E-04      & 1.34E-03 \\
1.00E-03      & 2.94E-06      & 1.40E-05 \\
\hline
\end{tabular}
\end{center}
\end{table}

\subsection{Equations Describing the Motion of the Air}
As seen from the experimental evidence, the velocities of air
encountered during coughing and sneezing never exceed a Mach-number
of $Ma=0.1$. Therefore, the air may be assumed as a Newtonian,
incompressible liquid, where buoyancy effects are modeled via
the Boussinesq approximation. The equations describing the
conservation of momentum, mass and energy for incompressible, 
Newtonian flows may be written as

$$ \rho \vvec_{,t} + \rho \vvec \cdot \Grd \vvec + \Grd p = 
 \nabla \cdot \mu \nabla \vvec + \rho \gvec + \beta \rho \gvec(T - T_0) 
      + \svec_v   ~~, 
                                                         \eqno(4.1.1) $$
$$                             \Div \vvec = 0      ~~,   \eqno(4.1.2) $$
$$ \rho c_p T_{,t} + \rho c_p \vvec \cdot \Grd T = 
       \nabla \cdot k \nabla T + s_e ~~.
                                                         \eqno(4.1.3) $$

\noi
Here $\rho, \vvec, p, \mu, \gvec, \beta, T, T_0, c_p, k$ denote 
the density, velocity vector, pressure, viscosity, gravity vector, 
coefficient of thermal expansion, temperature, reference temperature,
specific heat coefficient and conductivity respectively, and
$\svec_v, s_e$ momentum and energy source terms (e.g. due to particles
or external forces/heat sources).
For turbulent flows both the viscosity and the conductivity are
obtained either from additional equations or directly via a
large eddy simulation (LES) assumption through monotonicity induced
LES (MILES) \cite{Bor92,Fur99,Gri02,Ide19,Gim21,Ide21}. \\
The pathogen concentration is given by an advection-diffusion equation 
of the form:

$$ c_{,t} + \vvec \cdot \Grd c = 
       \nabla \cdot k_c \nabla c + s_c ~~,
                                                         \eqno(4.1.4) $$

\noi
where $c, k_c$ denote the concentration (pathogens/volume) and 
diffusivity of the pathogen, and $s_c$ is the source (or sink) term
(due to exhalation or inhalation).
In addition, a series of additional `diagnostics' equations may be
required. One of them is the `age of air' $\tau$ (a good measure for 
ventilation efficiency), given by:

$$ \tau_{,t} + \vvec \cdot \Grd \tau = 1 ~~.
                                                         \eqno(4.1.5) $$

\subsection{Numerical Integration of the Motion of the Air}
\label{numintnavto}
The last six decades have seen a large number of schemes that may be
used to solve numerically the incompressible Navier-Stokes
equations given by Eqns.(4.1.1-4.1.3). In the present case, the 
following design criteria were implemented:
\begin{itemize}
\item[-] Spatial discretization using {\bf unstructured grids} 
(in order to allow for arbitrary geometries and adaptive refinement);
\item[-] Spatial approximation of unknowns with 
{\bf simple linear finite elements} (in order to have a simple 
input/output and code structure);
\item[-] Edge-based data structures (for reduced access to memory and
indirect addressing);
\item[-] Temporal approximation using {\bf implicit integration of viscous
terms and pressure} (the interesting scales are the ones associated with
advection);
\item[-] Temporal approximation using {\bf explicit, high-order 
integration of advective terms};
\item[-] {\bf Low-storage, iterative solvers} for the resulting systems of
equations (in order to solve large 3-D problems); and
\item[-] Steady results that are {\bf independent from the timestep} chosen
(in order to have confidence in convergence studies).
\end{itemize} 
\noindent
The resulting discretization in time is given by the following projection
scheme \cite{Loh04,Loh06,Loh08}:
\begin{itemize}
\item[-] \ub{Advective-Diffusive Prediction}:
$\vvec^n, p^n \rightarrow \vvec^{*}$

$$ \svec' = - \Grd p^n + \rho \gvec 
          + \beta \rho \gvec (T^n - T_0) + \svec_v ~~,
\eqno(4.2.1)
$$

$$
\rho \vvec^i = \rho \vvec^n + \alpha^i \gamma \dt \left(
 - \rho \vvec^{i-1} \cdot \Grd \vvec^{i-1} 
 + \nabla \cdot \mu \nabla \vvec^{i-1} + \svec' \right)  ~~; ~~i=1,k-1~~;
\eqno(4.2.2)
$$

$$
 \left[ { \rho \over \dt} - \theta \nabla \cdot \mu \nabla \right]
   \left( \vvec^{k} - \vvec^n \right)
 + \rho \vvec^{k-1} \cdot \Grd \vvec^{k-1} = 
   \nabla \cdot \mu \nabla \vvec^{k-1} + \svec' ~~.  \eqno(4.2.3)
$$

\ms \noi
\item[-] \ub{Pressure Correction}: $p^n \rightarrow p^{n+1}$

$$
 \Div \vvec^{n+1} = 0                       ~~; \eqno(4.2.4)
$$
$$
 \rho {{ \vvec^{n+1} - \vvec^{*} }\over \dt} + \Grd ( p^{n+1} - p^n )
   = 0                                      ~~; \eqno(4.2.5)
$$

\noi
\item[ ] which results in

$$
 \Div { \dt \over \rho} \nabla ( p^{n+1} - p^n ) = \Div \vvec^{*} ~~;
\eqno(4.2.6)
$$

\ms \noi
\item[-] \ub{Velocity Correction}:
$\vvec^{*} \rightarrow \vvec^{n+1}$

$$
 \vvec^{n+1} = \vvec^{*} - { \dt \over \rho} \Grd ( p^{n+1} - p^n ) ~~.
\eqno(4.2.7)
$$
\end{itemize}

\noi
$\theta$ denotes the implicitness-factor for the viscous
terms ($\theta=1$: 1st order, fully implicit, $\theta=0.5$: 2nd order,
Crank-Nicholson). 
$\alpha^i$ are the standard low-storage Runge-Kutta coefficients
$\alpha^i=1/(k+1-i)$. The $k-1$ stages of Eqn.(4.2.2) may be seen as a 
predictor (or replacement)
of $\vvec^n$ by $\vvec^{k-1}$. The original right-hand side has not been
modified, so that at steady-state $\vvec^n=\vvec^{k-1}$, preserving the
requirement that the steady-state be independent of the timestep $\dt$.
The factor $\gamma$ denotes the local ratio of the stability limit for
explicit timestepping for the viscous terms versus the timestep chosen.
Given that the advective and viscous timestep limits are proportional to:

$$ \dt_a \approx {h \over {|\vvec|}} ~~;~~
   \dt_v \approx {{\rho h^2} \over \mu} ~~, \eqno(4.2.8)
$$

\noi
we immediately obtain

$$ \gamma = {{\dt_v} \over {\dt_a}}
    \approx {{\rho |\vvec| h }\over{\mu}} \approx Re_h  ~~,
\eqno(4.2.9)
$$

\noi
or, in its final form:

$$ \gamma = min(1,Re_h) ~~. \eqno(4.2.10) $$

\noi
In regions away from boundary layers, this factor is $O(1)$, implying 
that a high-order Runge-Kutta scheme is recovered. Conversely, for 
regions where $Re_h=O(0)$, the scheme reverts back to the usual
1-stage Crank-Nicholson scheme. 
Besides higher accuracy, an important benefit of explicit multistage 
advection schemes is the larger timestep one can employ. The increase in 
allowable timestep is roughly proportional to the number of stages used 
(and has been exploited extensively for compressible flow simulations 
\cite{Jam81}). 
Given that for an incompressible solver of the projection type 
given by Eqns.(4.2.1-4.2.7) most of the CPU time is spent solving the 
pressure-Poisson system Eqn.(4.2.6), the speedup
achieved is also roughly proportional to the number of stages used. \\
At steady state, $\vvec^{*}=\vvec^n=\vvec^{n+1}$ and the residuals of 
the pressure correction vanish,
implying that the result does not depend on the timestep $\dt$. \\
The spatial discretization of these equations is carried out via 
linear finite elements. The
resulting matrix system is re-written as an edge-based solver, allowing
the use of consistent numerical fluxes to stabilize the advection and
divergence operators \cite{Loh08}. \\
The energy (temperature) equation (Eqn.(4.1.3)) is integrated in a 
manner similar to the advective-diffusive prediction (Eqn.(4.2.2)), 
i.e. with an explicit, high order Runge-Kutta scheme for the advective 
parts and an implicit, 2nd order Crank-Nicholson scheme for the 
conductivity. 

\subsection{Immersed Body Techniques}
The information required from computational crowd dynamics (CCD) 
codes consists of
the pedestrians in the flowfield, i.e. their position, velocity,
temperature, as well inhalation and exhalation. As the CCD 
codes describe the pedestrians as points,
circles or ellipses, a way has to be found to transform this data
into 3-D objects. Two possibilities have been pursued here:
\begin{itemize}
\item{a)} Transform each pedestrian into a set of (overlapping) spheres that
approximate the body with maximum fidelity with the minimum amount of
spheres;
\item{b)} Transform each pedestrian into a set of tetrahedra that
approximate the body with maximum fidelity with the minimum amount of
tetrahedra.
\end{itemize}
The reason for choosing spheres or tetrahedra is that due to their
geometric simplicity one can perform
the required interpolation/ information transfer much faster than with
other polyhedra or geometric shapes. \\
In order to `impose' on the flow the presence of a pedestrian the
immersed boundary methodology is used. The key idea is to prescribe at
every CFD point covered by a pedestrian the velocity and temperature
of the pedestrian. For the CFD code, this translates into an extra
set of boundary conditions that vary in time and space as the
pedestrians move. This is by now a mature technology (see, e.g.
chapter 18 in \cite{Loh08} and the references cited therein). Fast search
techniques as well as extensions to higher order boundary conditions
may be found in \cite{Loh08,Loh08b}. Nevertheless, as the pedestrians
potentially change location at every timestep, the search for and the
imposition of new boundary conditions can add a considerable amount
of CPU as compared to `flow-only' runs.

\section{Modeling of Pedestrian Motion}
\label{modpedmotion}
The modeling of pedestrian motion has been the focus of research and
development for more than two decades \cite{Fru71,Pre71,Pel08}.
If one is only interested in
average quantities (average density, velocity), continuum models
\cite{Hug03} are an option. For problems requiring more realism,
approaches that model each individual are required \cite{Tha07}.
Among these, discrete space models (such as cellular automata
\cite{Blu98,Blu02,Tek00,Dij02,Scha02,Kes02,Klu03,Cou05,Lan06}), 
force-based models (such as
the social force model \cite{Hel95,Hel02,Qui03,Lak05,Loh10})
and agent-based techniques 
\cite{Pel06,Sud07,Guy09,Guy10,Vig10,Tor12,Cur12} have been 
explored extensively.
Together with insights from psychology and neuroscience (e.g.
\cite{Vis95,Tor12}) it has become clear that any
pedestrian motion algorithm that attempts to model reality should be
able to mirror the following empirically known facts and behaviours:
\begin{itemize}
\item[-] Newton's laws of motion apply to humans as well:
from one instant to another, we can only move within certain bounds
of acceleration, velocity and space;
\item[-] Contact between individuals occurs for high densities; these
forces have to be taken into account;
\item[-] Humans have a mental map and plan on how they desire to move
globally (e.g. first go here, then there, etc.); 
\item[-] Human motion is therefore governed by strategic 
(long term, long distance), tactical (medium
term, medium distance) and operational (immediate) decisions;
\item[-] In even moderately crowded situations of one person per
square meter (i.e. $O(1~p/m^2)$),
humans have a visual horizon of $O(2.5-5.0 m)$, and a perception range
of 120 degrees; thus, the influence of other humans beyond these
thresholds is minimal;
\item[-] Humans have a `personal comfort zone'; it is dependent on
culture and varies from individual to individual, but it cannot
be ignored;
\item[-] Humans walk comfortably at roughly 2 paces per second
(frequency: $\nu=2~Hz$); they are able to change the frequency for short
periods of time, but will return to $2~Hz$ whenever possible.
\end{itemize}
We remark that many of the important and groundbreaking work cited
previously took place within the gaming/visualization
community, where the emphasis is on `looking right'. Here, the aim
is to answer civil engineering or safety questions such as maximum
capacity, egress times under emergency, or comfort. Therefore,
comparisons with experiments and actual data are seen as essential
\cite{Loh10,Ise14a,Ise14b}.

\subsection{The PEDFLOW Model}
The PEDFLOW model \cite{Loh10} incorporates these requirements
as follows: individuals move according to Newton's laws of motion;
they follow (via will forces) `global movement targets'; at the local
movement level, the motion also considers the presence of other
individuals or obstacles via avoidance forces (also a type of will
force) and, if applicable, contact forces.
Newton's laws:

$$ m {{d \vvec}\over {dt}} = \fvec ~~,~~ 
     {{d \xvec}\over {dt}} = \vvec ~~, \eqno(5.1.1) $$

\noi
where $m, \vvec, \xvec, \fvec, t$ denote, respectively, mass,
velocity, position, force and time, are integrated in time
using a 2nd order explicit timestepping technique.
The main modeling effort is centered on $\fvec$.
In the present case the forces are separated into internal
(or will) forces [I would like to move here or there] and
external forces [I collided with another pedestrian or an
obstacle]. For the sake of completeness, we briefly
review the main forces used. For more
information, as well as verification and validation studies, see
\cite{Loh10,Ise14a,Ise14b,Zha14,Ise16,Ise16a,Ise16b,Loh16}.

\subsubsection{Will Force}
Given a desired velocity $\vvec_d$ and
the current velocity $\vvec$, this force will be of the form

$$ \fvec_{will} = g_w \left( \vvec_d - \vvec \right) ~~. 
                                                    \eqno(5.1.1.1) $$

\noi
The modelling aspect is included in the function $g_w$, which, in the
non-linear case, may itself be a function of $\vvec_d - \vvec$.
Suppose $g_w$ is constant, and that only the will force is acting.
Furthermore, consider a pedestrian at rest. In this case, we have:

$$  m {{d \vvec} \over {dt}} = g_w \left( \vvec_d - \vvec \right) 
                        ~~,~~ \vvec(0)=0 ~~,        \eqno(5.1.1.2) $$

\noi
which implies:

$$ \vvec = \vvec_d \left( 1 - e^{-\alpha t} \right) ~~,~~
  \alpha = {{g_w}\over{m}} = {1 \over t_r} ~~, \eqno(5.1.1.3) $$

\noi
and

$$ {{d \vvec} \over {dt}}(t=0) = \alpha \vvec_d 
                               = {\vvec_d \over t_r} ~~. 
                                                \eqno(5.1.1.4) $$

\noi
One can see that the crucial parameter here is the `relaxation time'
$t_r$ which governs the initial acceleration and `time to desired
velocity'. Typical values are $\vvec_d=1.35~m/sec$ and $t_r=O(0.5~sec)$.
The `relaxation time' $t_r$ is clearly dependent on the fitness of the
individual, the current state of stress, desire to reach a
goal, climate, signals, noise, etc. Slim, strong individuals will have
low values for $t_r$, whereas obese or weak individuals will have high
values for $t_r$.
Furthermore, dividing by the mass of the individual allows all other
forces (obstacle and pedestrian collision avoidance, contact, etc.) to
be scaled by the `relaxation time' as well, simplifying the modeling
effort considerably.
\par \noi
The direction of the desired velocity

$$ \svec = {{\vvec_d} \over {|\vvec_d|}}         \eqno(5.1.1.5) $$

\noi
will depend on the type of pedestrian and the cases considered. A
single individual will have as its goal a desired position
$\xvec_d(t_d)$ that he/she would like to reach at a certain time $t_d$.
If there are no time constraints, $t_d$ is simply set to a large number.
Given the current position $\xvec$, the direction of the velocity is
given by

$$ \svec = {{ \xvec_d(t_d) - \xvec } \over {|\xvec_d(t_d) - \xvec|} } ~~,
                                                    \eqno(5.1.1.6) $$

\noi
where $\xvec_d(t_d)$ denotes the desired position (location, goal) of
the pedestrian at the desired time of arrival $t_d$. 
For members of groups, the goal is always to stay close to the leader.
Thus, $\xvec_g(t_g)$ becomes the position of the leader. In the case of
an evacuation simulation, the direction is given by the gradient of the
perceived time to exit $\tau_e$ to the closest perceived exit:

$$ \svec = {{ \nabla \tau_e } \over {|\nabla \tau_e}| } ~~. 
                                                    \eqno(5.1.1.7) $$

\noi
The magnitude of the desired velocity $|\vvec_d|$ depends on the fitness
of the individual, and the motivation/urgency to reach a certain place
at a certain time. Pedestrians typically stroll leisurely at
$0.6-0.8~m/sec$, walk at $0.8-1.0~m/sec$, jog at $1.0-3.0~m/sec$,
and run at $3.0-10.0~m/sec$.

\subsubsection{Pedestrian Avoidance Forces}
The desire to avoid collisions with other individuals is modeled by
first checking if a collision will occur. If so, forces are applied
in the direction normal and tangential to the intended motion. The
forces are of the form:

$$ f_i = f_{max}/( 1 + \rho^p) ~~;~~ 
   \rho  = | \xvec_i - \xvec_j |  /  r_i ~~, \eqno(5.1.2.1) $$

\noi
where $x_i$, $x_j$ denote the positions of individuals $i,j$,
$r_i$ the radius of individual $i$, and $f_{max}=O(4)f_{max}(will)$.
Note that the forces weaken with increasing non-dimensional distance
$\rho$. For years we have used $p=2$, but, obviously, this can (and
probably will) be a matter of debate and speculation (perhaps a
future experimental campaign will settle this issue).
In the far range, the forces are mainly orthogonal to the direction of
intended motion: humans tend to move slightly sideways without decelerating.
In the close range, the forces are also in the direction of intended
motion, in order to model the slowdown required to avoid a collision.

\subsubsection{Wall Avoidance Forces}
Any pedestrian modeling software requires a way to input geographical
information such as walls, entrances, stairs, escalators, etc. In
the present case, this is accomplished via a triangulation (the so-called
background mesh). A distance to walls map (i.e. a function $d_w(x)$
is constructed using fast marching techniques on unstructured grids),
and this allows to define a wall avoidance force as follows:

$$ \fvec = - f_{max}{ 1 \over { 1 + ({d_w \over r})^p }} 
         \cdot \nabla d_w ~~,~~p=2             \eqno(5.1.3.1) $$

\noi
Note that $|\nabla d_w|=1$. The default for the maximum wall 
avoidance force is $f_{max}=O(8)f_{max}(will)$. The desire to be 
far/close to a wall also depends on cultural background.

\subsubsection{Contact Forces}
When contact occurs, the forces can increase markedly. Unlike will
forces, contact forces are symmetric. Defining

$$ \rho_{ij} = |\xvec_i - \xvec_j|/(r_i + r_j) ~~,   \eqno(5.1.4.1) $$

\noi
these forces are modeled as follows:

$$ \rho_{ij} < 1:  f = - [ f_{max} /( 1 +  \rho_{ij}^p)]~~; ~~  p=2
                                                \eqno(5.1.4.2a)  $$
$$ \rho_{ij} > 1:  f = - [2f_{max} /( 1 +  \rho_{ij}^p)]~~; ~~  p=2 
                                                \eqno(5.1.4.2b) $$

\noi
and $f_{max}=O(8)f_{max}(will)$.

\subsubsection{Motion Inhibition}
A key requirement for humans to move is the ability to put one foot
in front of the other. This requires space. Given the comfortable
walking frequency of $\nu=2~Hz$, one is able to limit the comfortable
walking velocity by computing the distance to nearest neighbors and
seeing which one of these is the most `inhibiting'.

\subsubsection{Psychological Factors}
The present pedestrian motion model also incorporates a number of 
psychological factors that,
among the many tried over the years, have emerged as important for
realistic simulations. Among these, we mention:
\begin{itemize}
\item[-] Determination/Pushiness: it is an everyday experience that in
crowds, some people exhibit a more polite behavior than others. This
is modeled in PEDFLOW by reducing the collision avoidance forces of
more determined or `pushier' individuals. Defining a determination or
pushiness parameter $\ptilde$, the avoidance forces are reduced by 
$(1-\ptilde)$. Usual ranges for $\ptilde$ are $0.2\le\ptilde\le0.8$.
\item[-] Comfort zone: in some cultures (northern Europeans are a
good example) pedestrians want to remain at some minimum distance from
contacting others. This comfort zone is an input parameter in PEDFLOW,
and is added to the radii of the pedestrians when computing collisions
avoidance and pre-contact forces.
\item[-] Right/Left Avoidance and Overtaking: in many western countries
pedestrians tend to avoid incoming pedestrians by stepping towards
their right, and overtake others on the left. However, this is not the
norm everywhere, and one has to account for it.
\end{itemize}

\subsection{Numerical Integration of the Motion of Pedestrians}
\label{numintpeds}
The equations describing the position and velocity of a pedestrian
may be formulated as a system of nonlinear Ordinary Differential 
Equations of the form:

$$ {{d\uvec_p} \over {dt}} = \rvec(\uvec_p, \xvec, \uvec_f) ~~,
                                                   \eqno(5.2.1)    $$

\noi
where $\uvec_p$ denote the variables of the pedestrians
(positions, velocities, ...) and $\rvec$ a right-hand-size that
depends on $\uvec_p$, the position of the pedestrian (e.g. geographical
obstacles) and the variables interpolated from the flow
domain (e.g. smoke) $\uvec_f$.
These ODEs are integrated with explicit Runge-Kutta schemes,
typically of order 2. \\
The geographic information required, such as terrain data
(inclination, soil/water, escalators, obstacles, etc.), climate data
(temperature, humidity, sun/rain, visibility), signs,
the location and accessibility of guidance personnel, as well as doors,
entrances and emergency exits is stored
in a so-called background grid consisting of triangular elements. This
background grid is used to define the geometry of the problem.
At every instance, a pedestrian will be located in one of the elements
of the background grid. Given this `host element' the geographic
data, stored at the nodes of the background grid, is interpolated
linearly to the pedestrian.
The closest distance to a wall $\delta_w$ or exit(s) for any given
point of the background grid evaluated via a fast ($O(N\ln(N))$)
nearest neighbour/heap list technique (\cite{Loh08,Loh10}).
For cases with visual or smoke impediments, the closest distance to
exit(s) is recomputed every few seconds of simulation time.

\subsection{Linkage to CFD Codes}
The information required from CFD codes such as 
FEFLO \cite{Loh01,Loh02,Loh13} consists of
the spatial distribution of temperature, smoke, other toxic or
movement impairing substances in space, as well as pathogen 
distribution. This information is
interpolated to the (topologically 2-D) background mesh at every
timestep in order to calculate properly the visibility/ reachability
of exits, routing possibilities, smoke, toxic substance or
pathogen inhalation, and any other flowfield variable required 
by the pedestrians. As the tetrahedral grid used for the CFD code 
and the triangular
background grid of the CCD code do not change in time, the interpolation
coefficients need to be computed just once at the beginning of the
coupled run. While the transfer of information from CFD to CCD is
voluminous, it is very fast, adding an insignificant amount
to the total run-times.

\section{Coupling Methodology}
\label{coupling}
\begin{figure}
	\includegraphics[width=10.0cm]{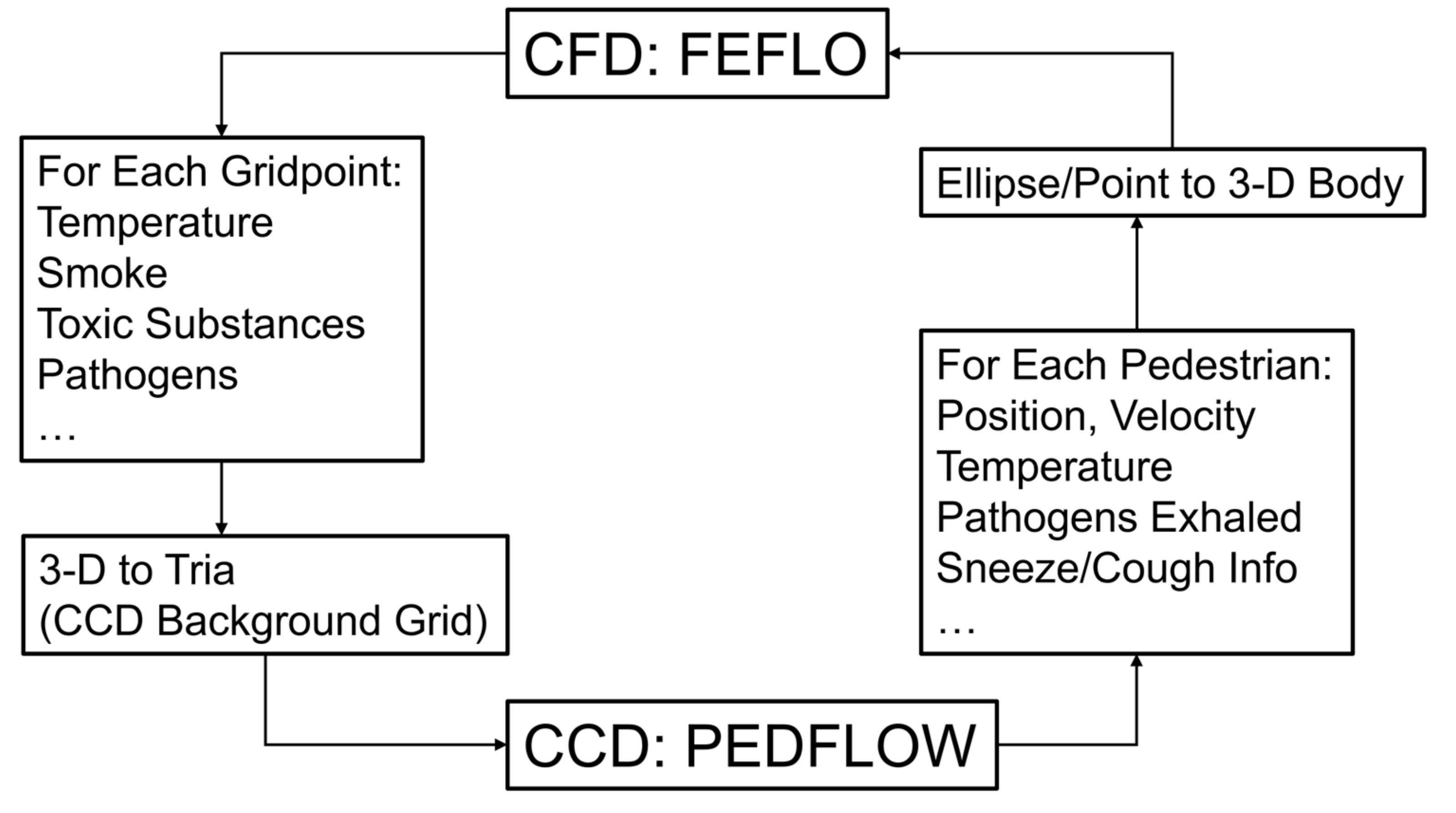}
	\caption{Coupling CFD and CCD Codes}
	\label{f:coupling}
\end{figure}

The coupling methodology used is shown in Figure~\ref{f:coupling}. 
The CFD code computes
the flowfield, providing such information as temperature, smoke, toxic
substance and pathogen concentration, and any other flow quantity 
that may affect the movement of pedestrians. These variables are then 
interpolated to the position where the pedestrians are, and are 
used with all other
pertinent information (e.g. will-forces, targets, exits, signs, etc.)
to update the position, velocity, inhalation of smoke, toxic substances 
or pathogens, state of exhaustion or intoxication, and any other
pertinent quantity that is evaluated for the pedestrians. The position,
velocity and temperature of the pedestrians, together with information
such as sneezing or exhaling air, is then transferred to the
CFD code and used to modify and update the boundary conditions of the
flowfield, particles and pathogen concentrations in the regions 
where pedestrians are present.

\noi
Of the many possible coupling options
(see e.g. \cite{Bun06,Sch08,Ceb05}),
we have implemented the simplest one: loose coupling with
sequential timestepping (\cite{Loh95,Loh11,Loh16a}). This is 
justified, as
the timesteps of both the flow and pedestrian solvers are very small,
so that possible coupling errors are negligible. PEDFLOW typically
runs with fixed timesteps of $\Delta t_{p}=0.05~sec$, while the 
timestep chosen by FEFLO depends on mesh size and velocity. Should
the timestep of FEFLO be less than the default value for PEDFLOW,
then PEDFLOW automatically reduces its timestep to be the same as
FEFLO. During the course of many cases run, we have never encountered 
any stability problems with this loose coupling and timestepping
strategy.

\subsection{Placement of Pathogen Loads in Space}
A background grid is used for the placement of geographical information
in PEDFLOW. The same grid can be used to track pathogen concentrations.
As infected pedestrians move through this grid, they exhale pathogen 
loads - either through sneezing, coughing, shouting or talking.
These pathogen loads are added to the concentration $c$ on the
background grid.

\subsection{Generation of Pathogen Loads}
Pathogen loads are generated whenever an infected pedestrian exhales,
either violently in bursts (e.g. sneezing, coughing, shouting), or
continuously (e.g. loud talking). The amount
of viral load can vary widely depending on the mode 
\cite{Lou67,Tan06,Xie07,Lin10,Teu10,Gup11b,Joh11,Tan11,Hal12,Lin12,Tan12,Mil13,Tan13,Wei16,Asa19,Asa20a,Asa20b,WHO20,To20}, 
the state of infection of the pedestrian, and many other factors (the 
term `superspreaders' has been used in the medical literature).
The position, velocity and orientation of pedestrians, as well
as their behaviour while sneezing or coughing is transmitted from
the pedestrian code to the flow code. The flow code then generates
the proper boundary conditions for the exhalation.
In order to simulate a sneeze/cough of duration $T_s$, the
velocity, temperature and pathogen concentration in a spherical 
region of radius ($r=5~cm$) near the pedestrian mouth is reset
at the beginning of each timestep according to the following
triangular function:

$$ f(t)= 
\begin{cases}
     {{2t} \over T_s} & if:         0 \le t \le   T_s/2 \cr
2 - {{2t}\over T_s} & if:   T_2/2 \le t \le T_s \cr
     0                  & if: T_s \le t              
\end{cases}
~~, \eqno(37) $$

$$ v(t) = v_{max} f(t)  ~~,~~ 
   T(t) = T_a + f(t) ( T_b - T_a ) ~~,~~
   c(t) = c_{max} f(t)  ~~,
    \eqno(38) $$

\noi
where $T_s=1$~[sec], $v_{max}=O(5)$~[m/sec], $T_a=20$~[$^o$C], 
$T_b=37$~[$^o$C], and $c_{max}$ an estimated concentration of virons 
exhaled that depends on the health of the pedestrian
(in the present case simply set to $c_{max}=1$ in arbitrary units).
The direction of the velocity is set from the orientation of the
pedestrian. The walking velocity of the pedestrian is then added
in order to obtain the final value. We remark that this simplified
model represents a cross-section of experimental data 
\cite{Cha09,Gup09,Gup10}, but that this is an active area of 
research \cite{Bus20}.

\subsection{Inhalation of Pathogen Loads}
As pedestrians walk or run through the clouds of viral loads, they
inhale a certain amount of viruses. Given the local concentration 
of viral load $c$ and the breathing rate of a pedestrian, the
total number of viruses inhaled can be integrated in time.
The assumption is made that once the inhaled viral load reaches
the infectious dose, the pedestrian is considered infected.

\section{Examples}
In the sequel, we show examples of different situations.
We remark that these are by no means
exhaustive or unique: the simulation of aerosol transmission via
high-fidelity CFD techniques has received considerable attention
in recent years, and has been carried out with commercial and
open source software worldwide (see, e.g. 
\cite{Ip07,Vil13,Zha17,Abu20,Dbo20a,Dbo20b,Li20,Loh20,Zoh20a,Loh21}).
The CFD code used is FEFLO, which was validated
for the class of problems considered here over many years
\cite{Ram96,Ram99,Cam04a,Cam06,Loh06,Loh14}.

\subsection{Fast Food Restaurant}
\begin{figure}
	\includegraphics[width=10.0cm]{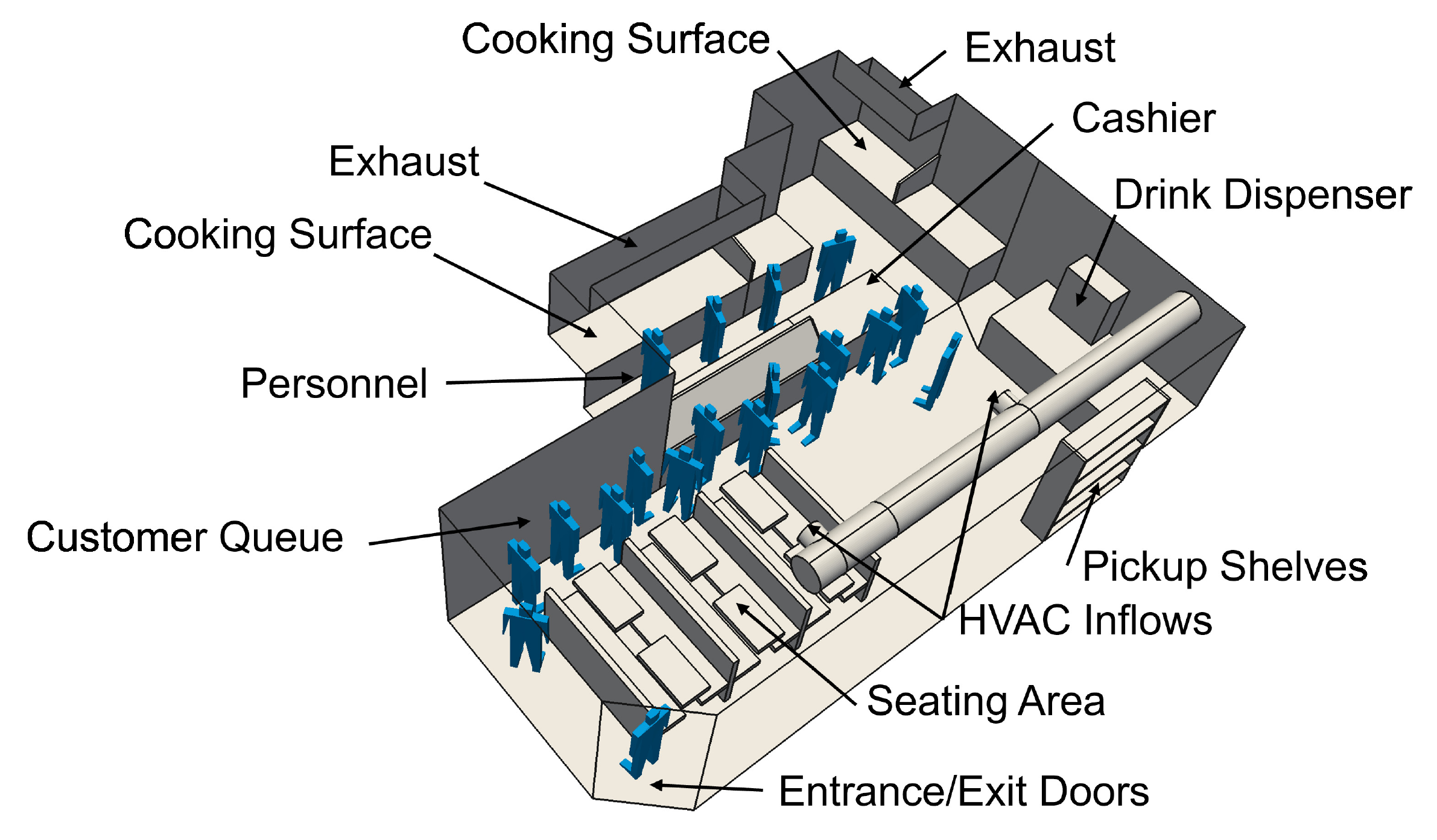}
	\caption{Fast Food Restaurant: Floorplan}
	\label{f:rest}
\end{figure}

This case considers a typical fast food restaurant.
The floorplan, geometry and general outlay are shown in Figure~\ref{f:rest}. 
Pedestrians enter the door at a rate of 0.1~p/sec, line up in the 
queue, order the food that can be seen behind the counter, and then 
pay at the cashier. The (random) `loiter times' at each of the 
stations along the counter with the food ranged from 
$t_{loit}=20-40~sec$.

\begin{figure}
	\includegraphics[width=6.0cm]{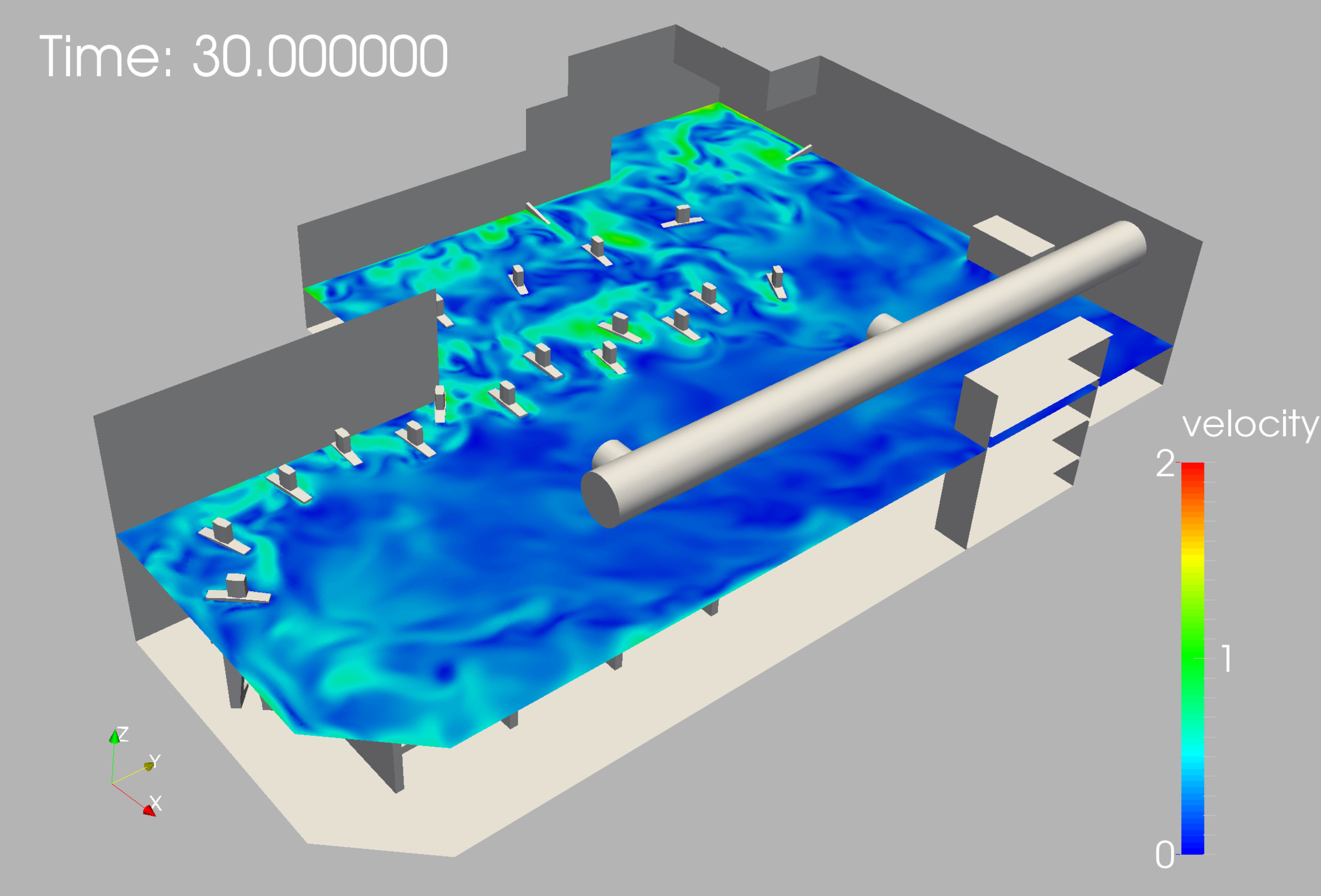}
	\includegraphics[width=6.0cm]{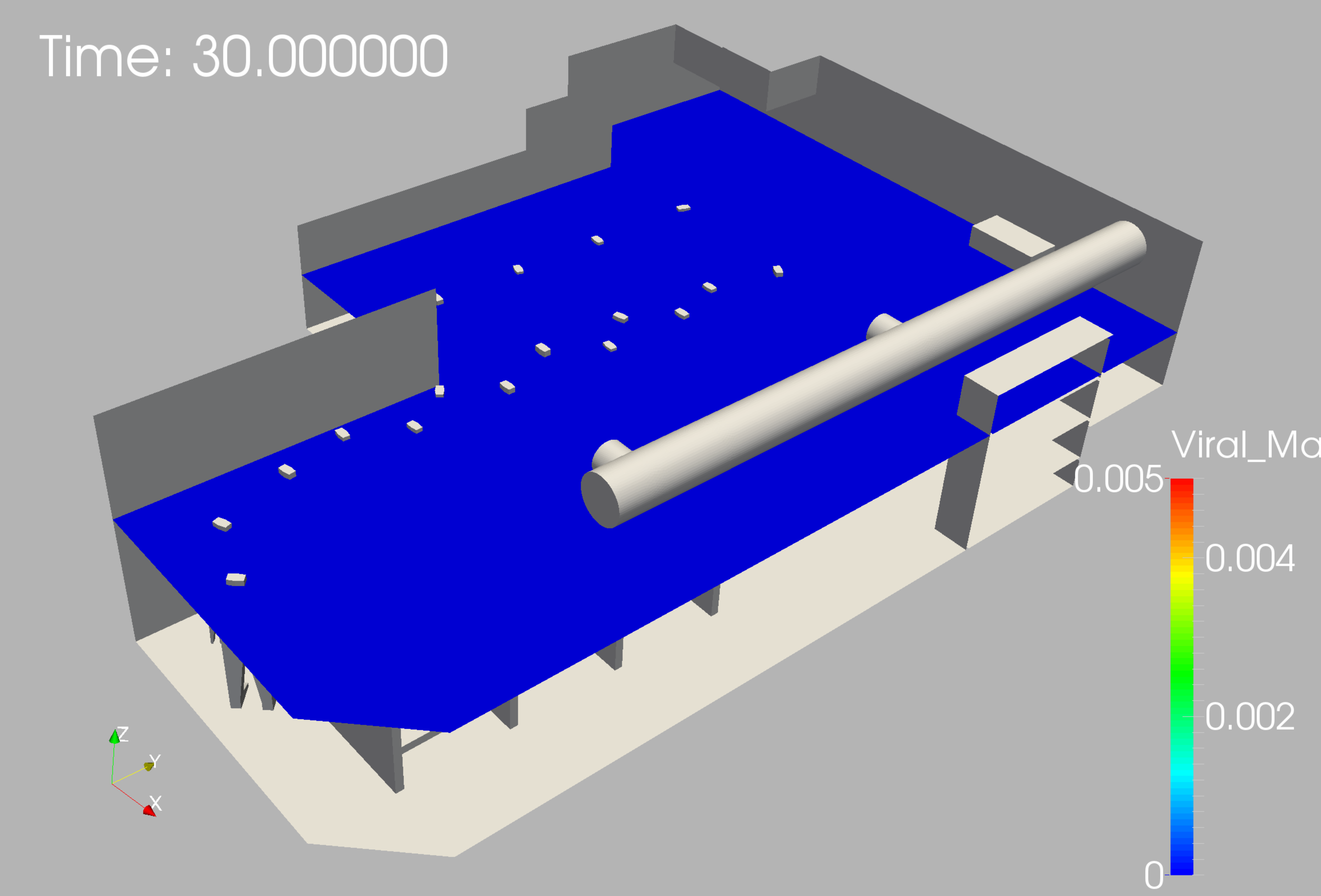}
	\caption{Fast Food Restaurant: Velocity and Viral Load
at $t=30.00~sec$.}
	\label{f:rest_1}
\end{figure}

\begin{figure}
	\includegraphics[width=6.0cm]{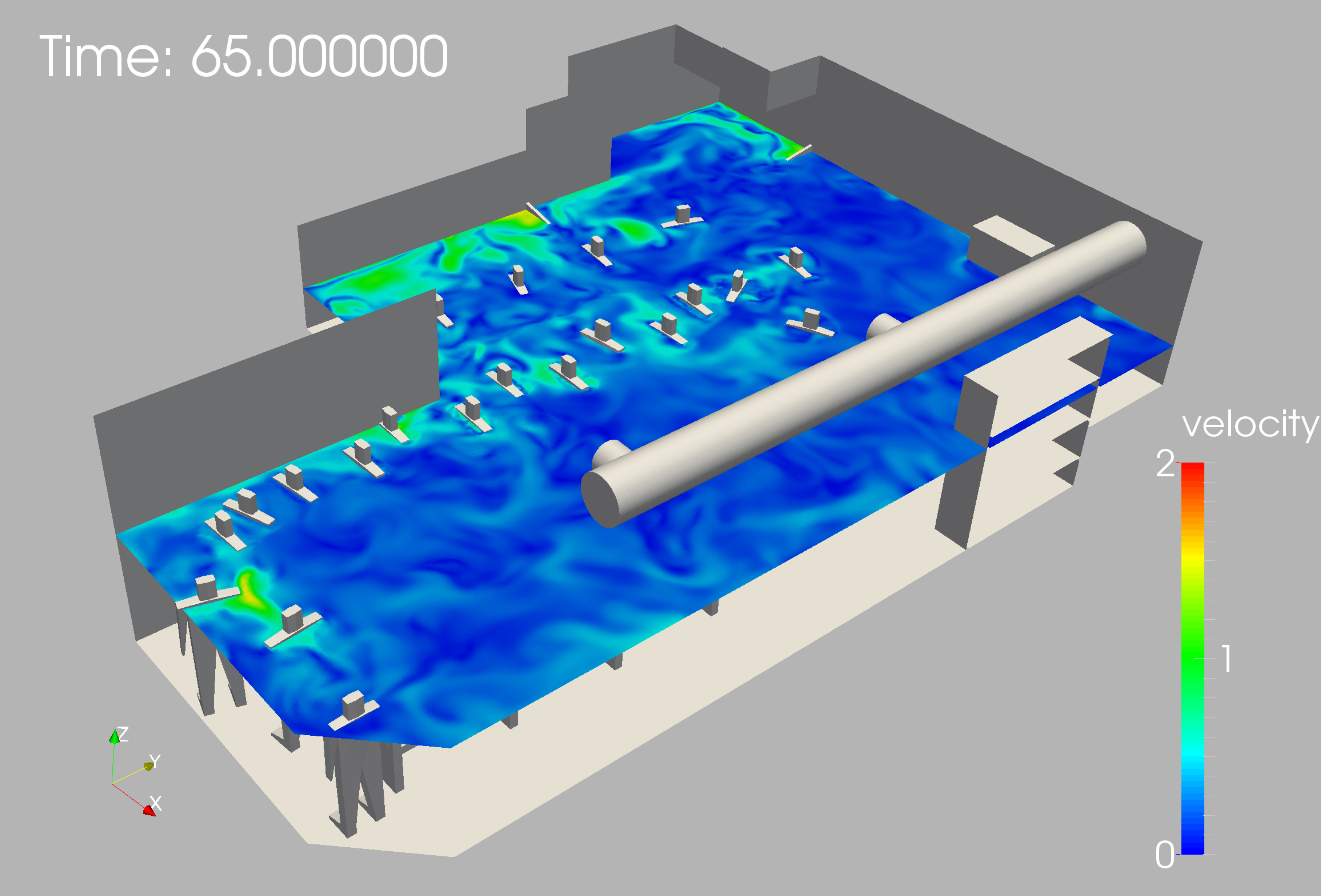}
	\includegraphics[width=6.0cm]{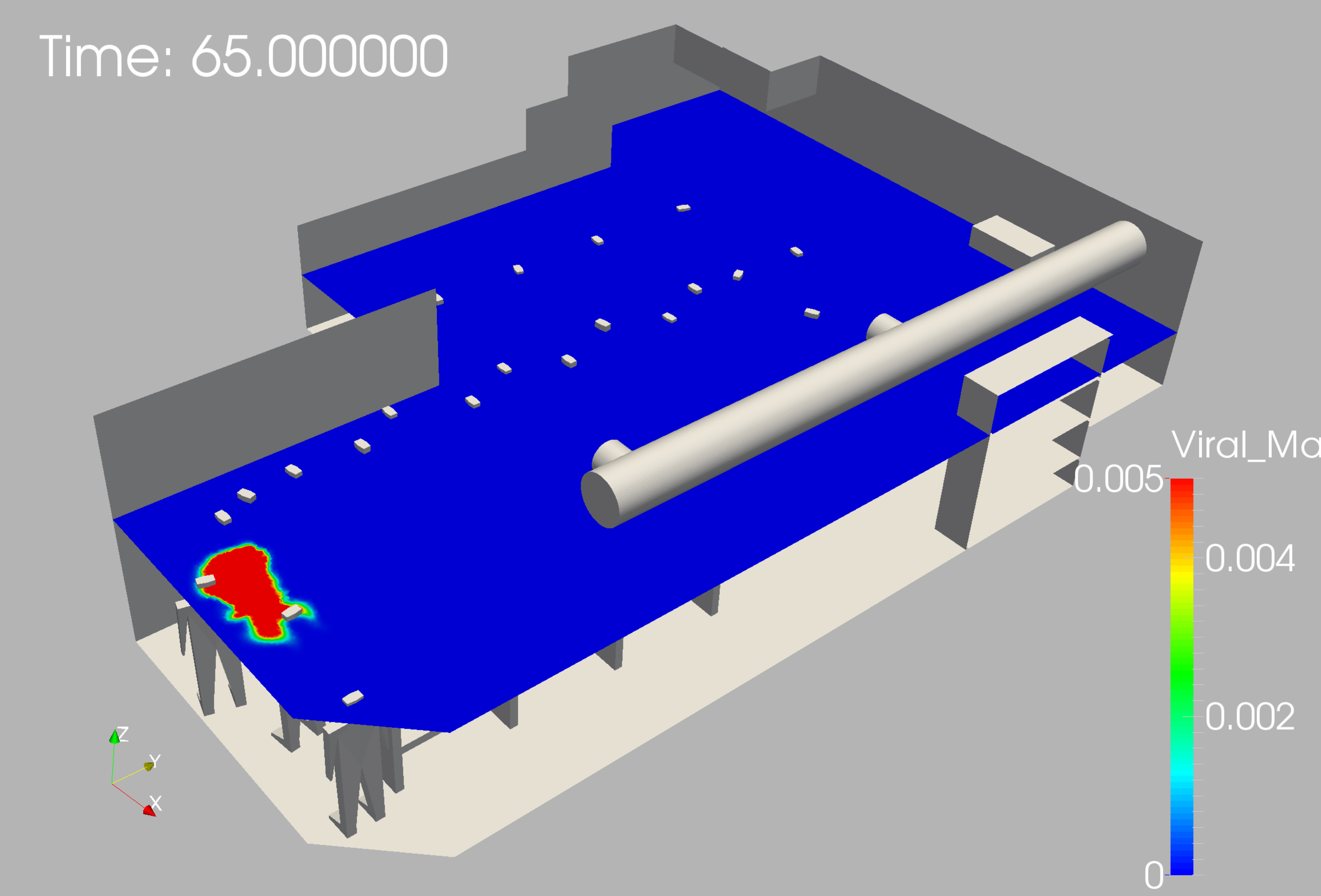}
	\caption{Fast Food Restaurant: Velocity and Viral Load
at $t=65.00~sec$}
	\label{f:rest_2}
\end{figure}

\begin{figure}
	\includegraphics[width=6.0cm]{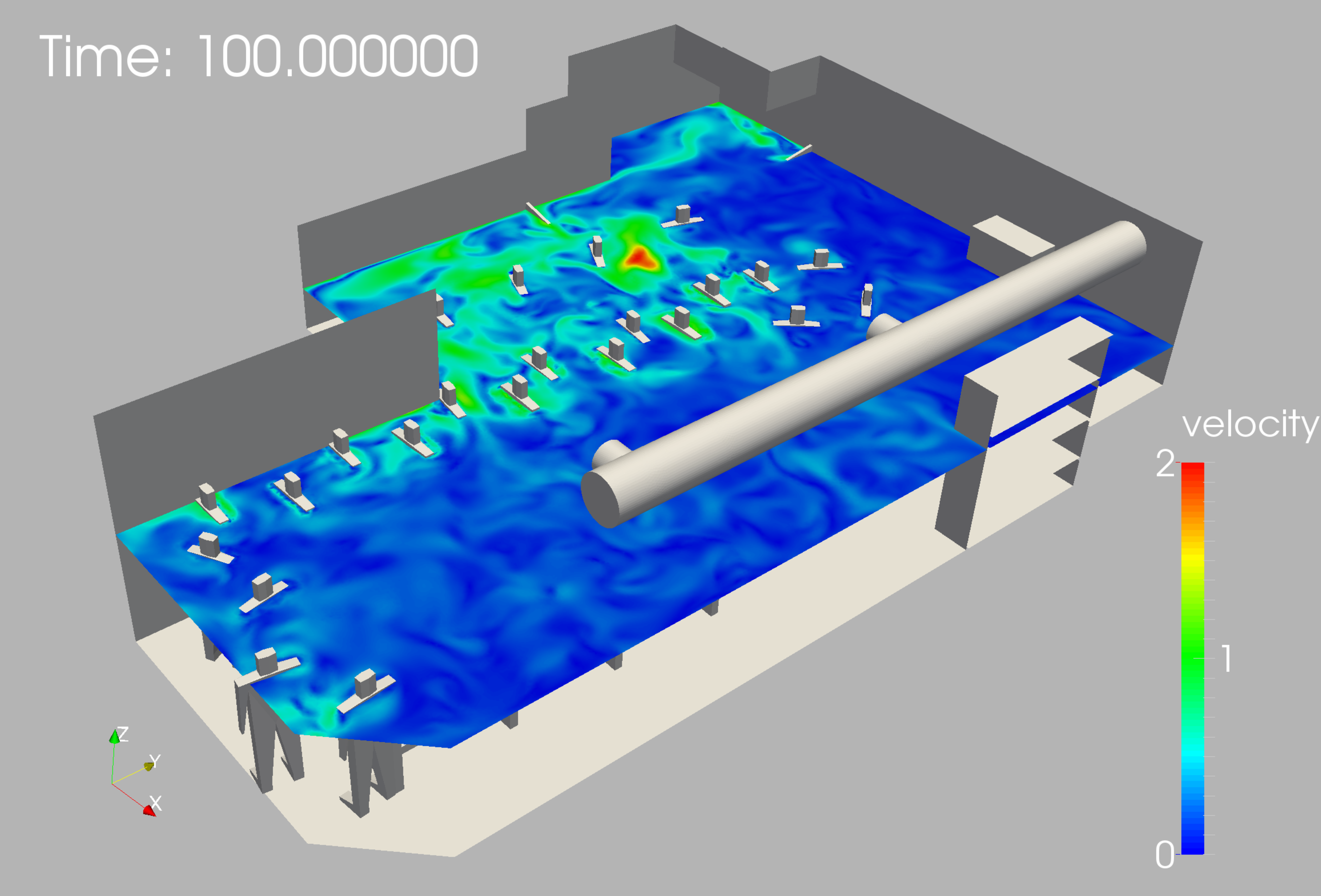}
	\includegraphics[width=6.0cm]{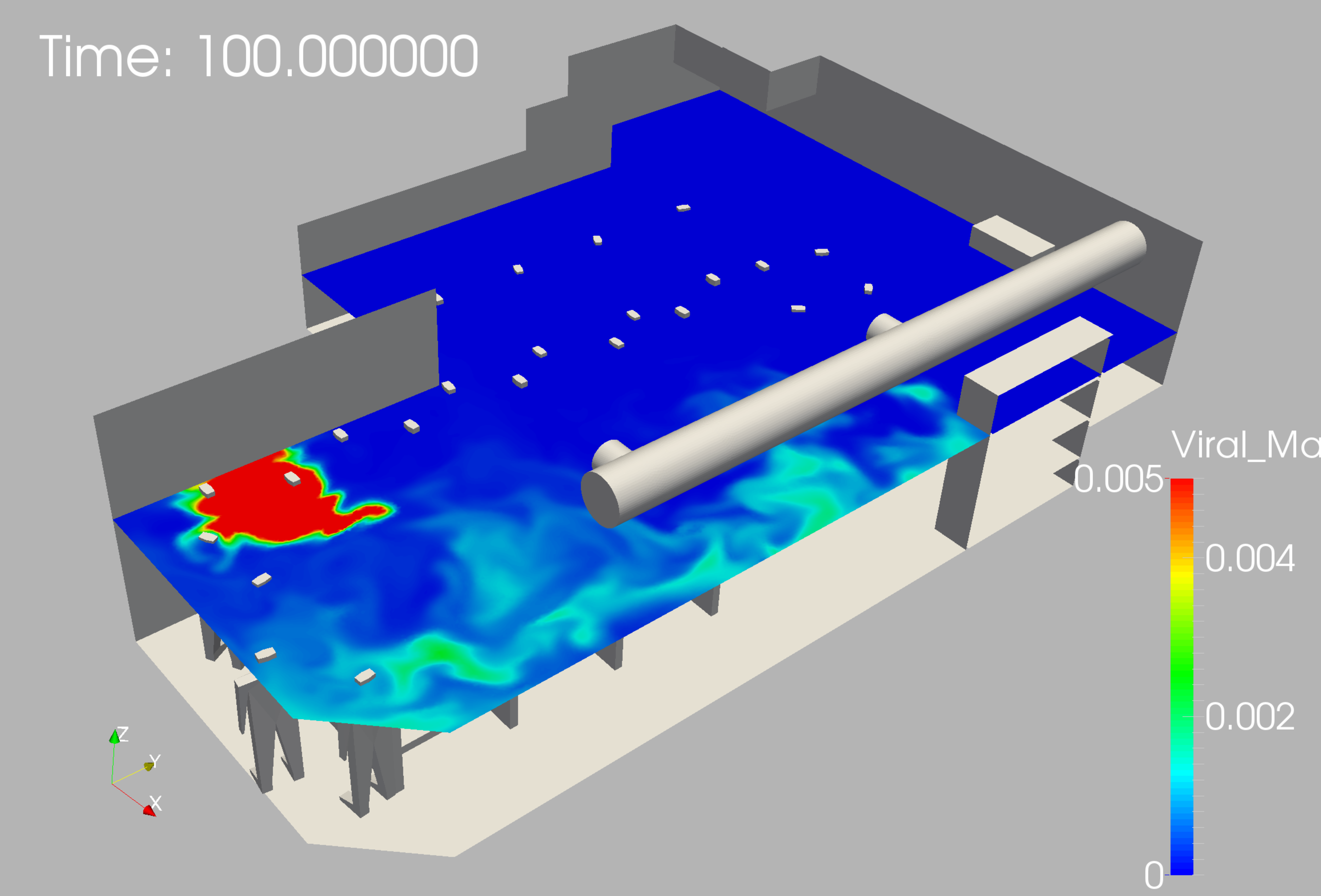}
	\caption{Fast Food Restaurant: Velocity and Viral Load
at $t=100.00~sec$}
	\label{f:rest_3}
\end{figure}

\ms \noi
The air in the room is refreshed via two HVAC inflow ducts close 
to the ceiling and two large exhaust surfaces located above the 
cooking surfaces. The initial airflow was considered at rest, and
the room temperature was assumed to be
$T_0=20^oC$. The air coming in from the two HVAC inflow ducts
has a velocity of $v_{HVAC}=2~m/sec$ and a temperature of
$T_{HVAC}=18^oC$. The heat released by the cooking surfaces was modeled
via a volumetric source that occupied 20~cm above the cooking surfaces
with a specific strength of $q=10^4~W/m^3$. At the beginning, all
pedestrians are considered healthy. Thereafter, 70\% of the pedestrians
entering are assigned as healthy, 10\% as infected but
asymptomatic (i.e. not infecting) and 20\% as infected and infecting.
For the infecting pedestrians, the average time between sneezing
is set to $\Delta t_s=30~sec$ with a variation of
$\delta t_s=5~sec$ and a sneeze duration of $T_s=1~sec$.
The mesh size in the region of the pedestrians was set to
$h=0.05~m$, which implies approximately 10 flow points per pedestrian
cross-section of $d_p=0.50~m$, sufficient in
order to obtain the vortical motion due to pedestrian movement or
flow obstruction.
The total mesh size was in excess of $1.8 \cdot 10^7$~elements.
The two minutes of simulated time took approximately 48~hours using
32~cores. The results obtained are shown in Figures~\ref{f:rest_1}--\ref{f:rest_3}.
Note the abscence of viral matter at the beginning. However, after
approximately two minutes (during which time two sneezing events
occured), a large area that is occupied
by pedestrians shows the presence of viral matter, highlighting the
need for better ventilation close to the entrance. The fluid dynamic
explanation is that the second ventilation outlet `blocks' the 
entrance region, producing a semi-stagnant region of recirculating 
flow where viral matter can build up.

\subsection{Long Passage in Train Station}
\begin{figure}
	\includegraphics[width=6.0cm]{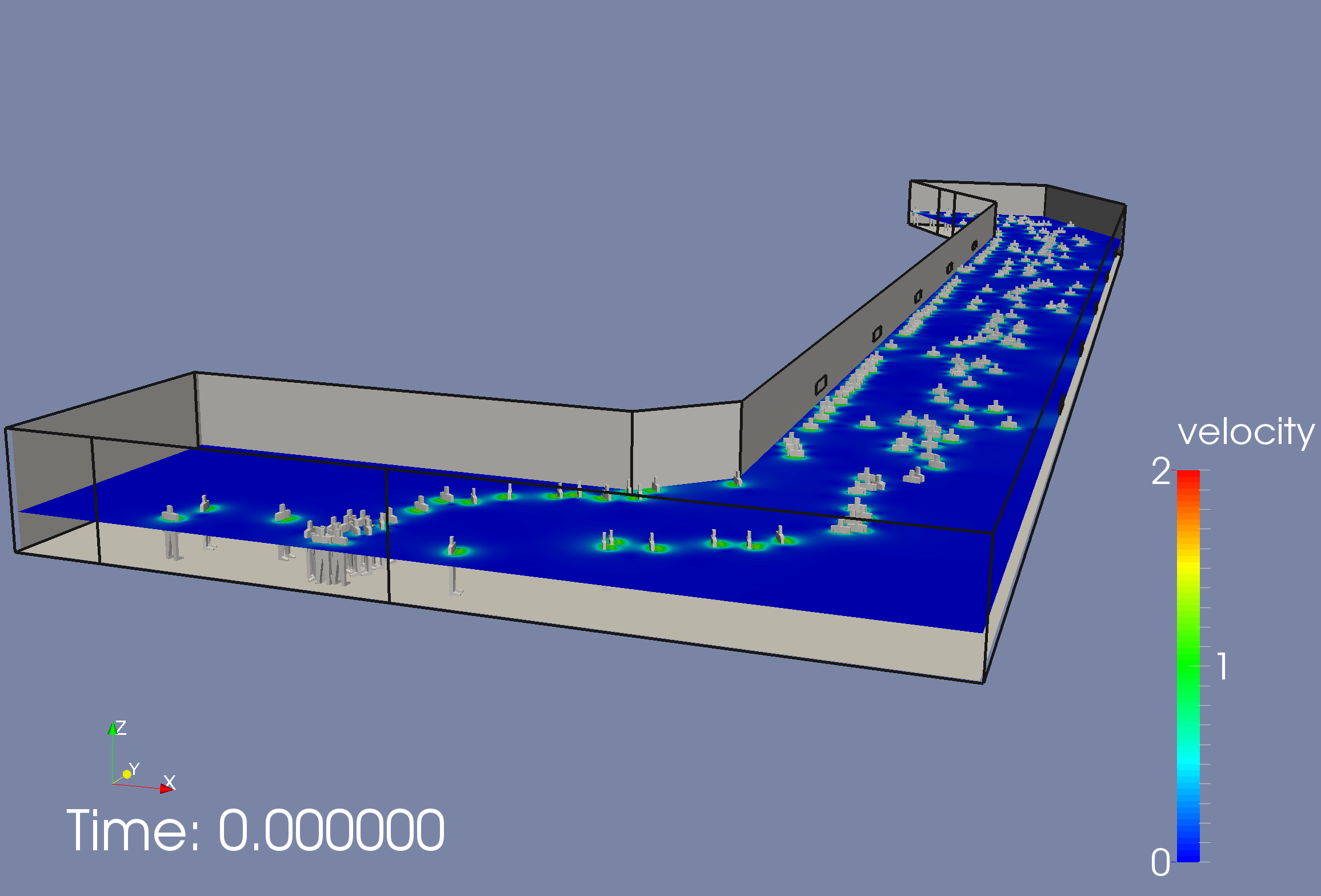}
	\includegraphics[width=6.0cm]{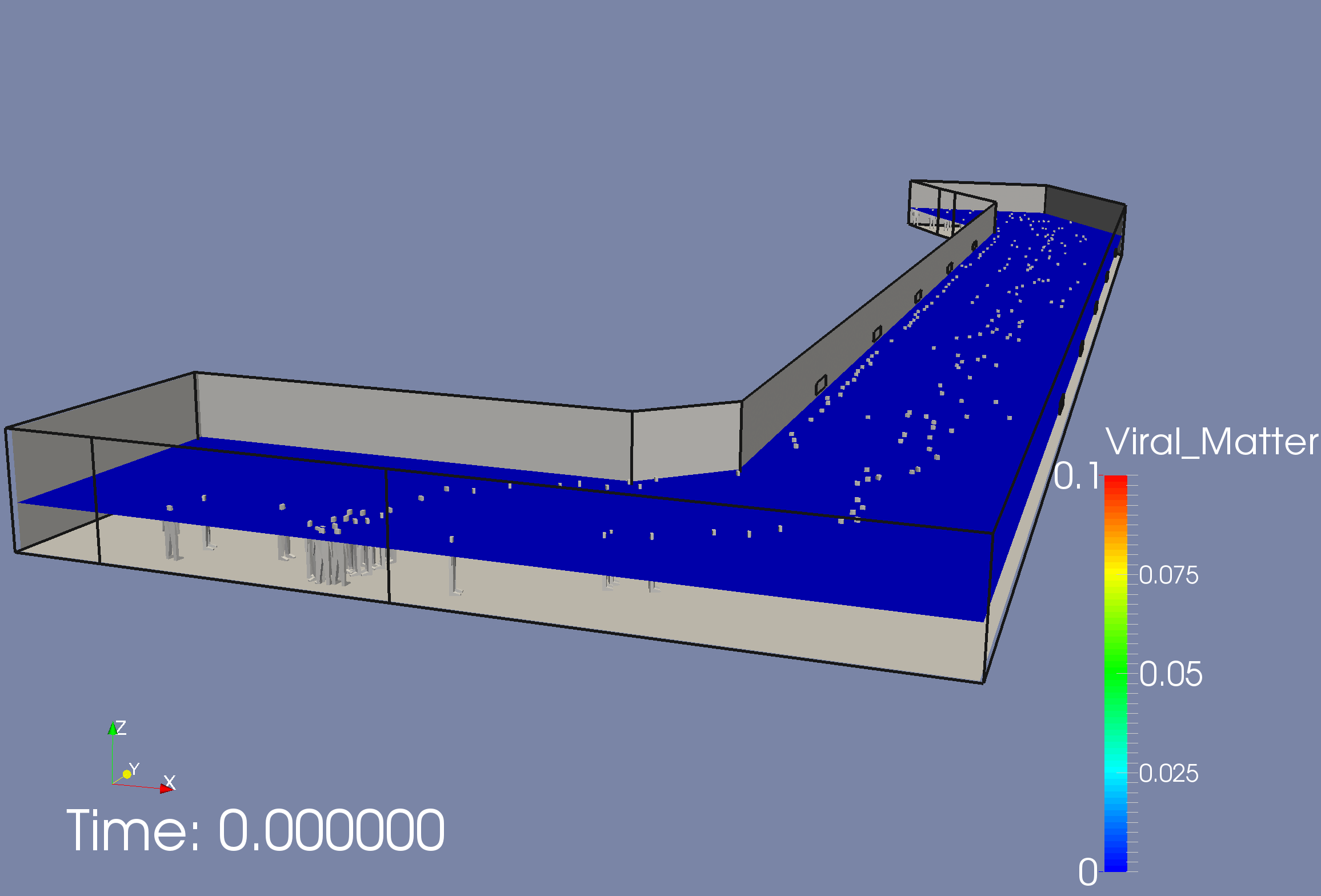}
	\caption{Long Passage: Initial Conditions}
	\label{f:pass_1}
\end{figure}

\begin{figure}
	\includegraphics[width=6.0cm]{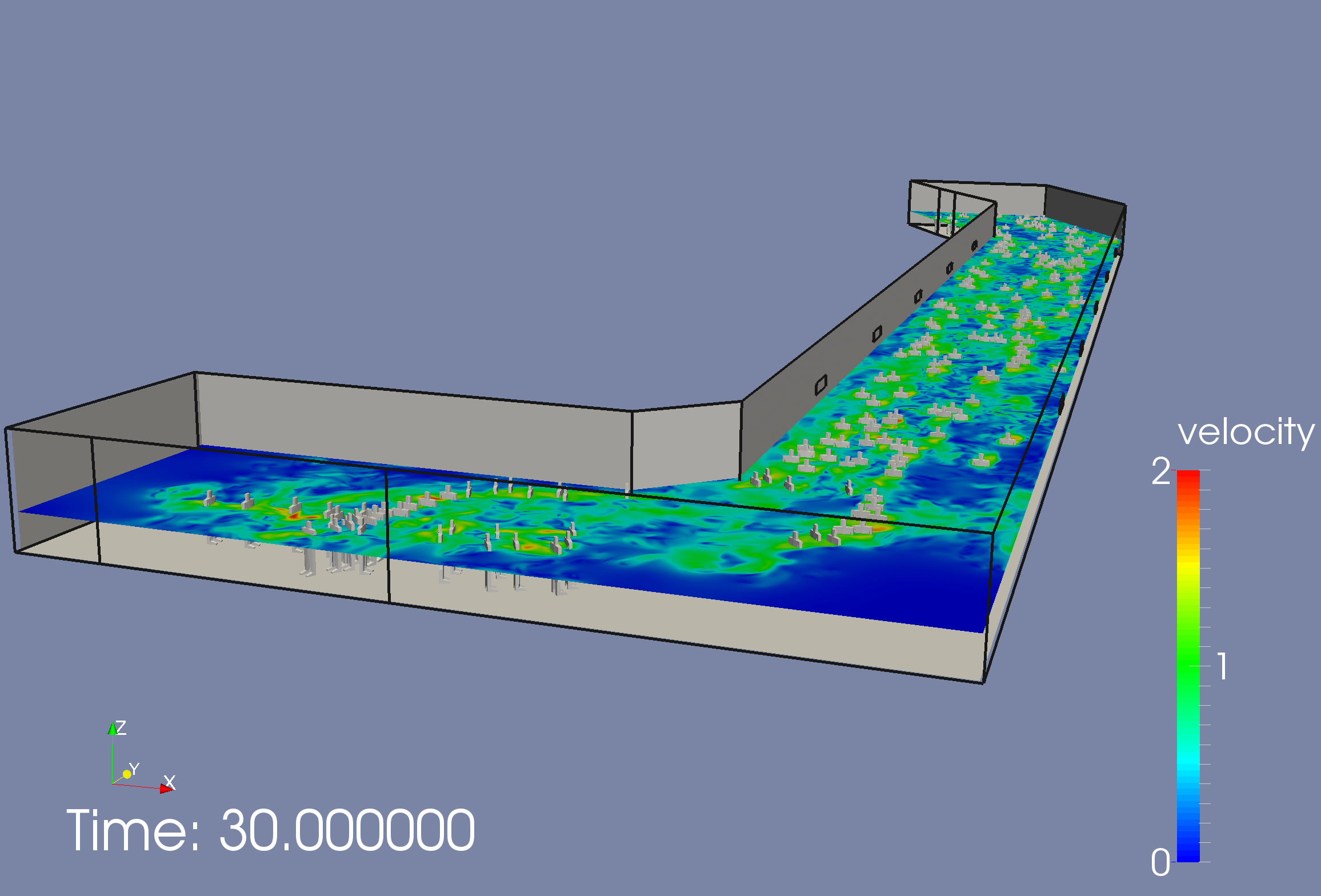}
	\includegraphics[width=6.0cm]{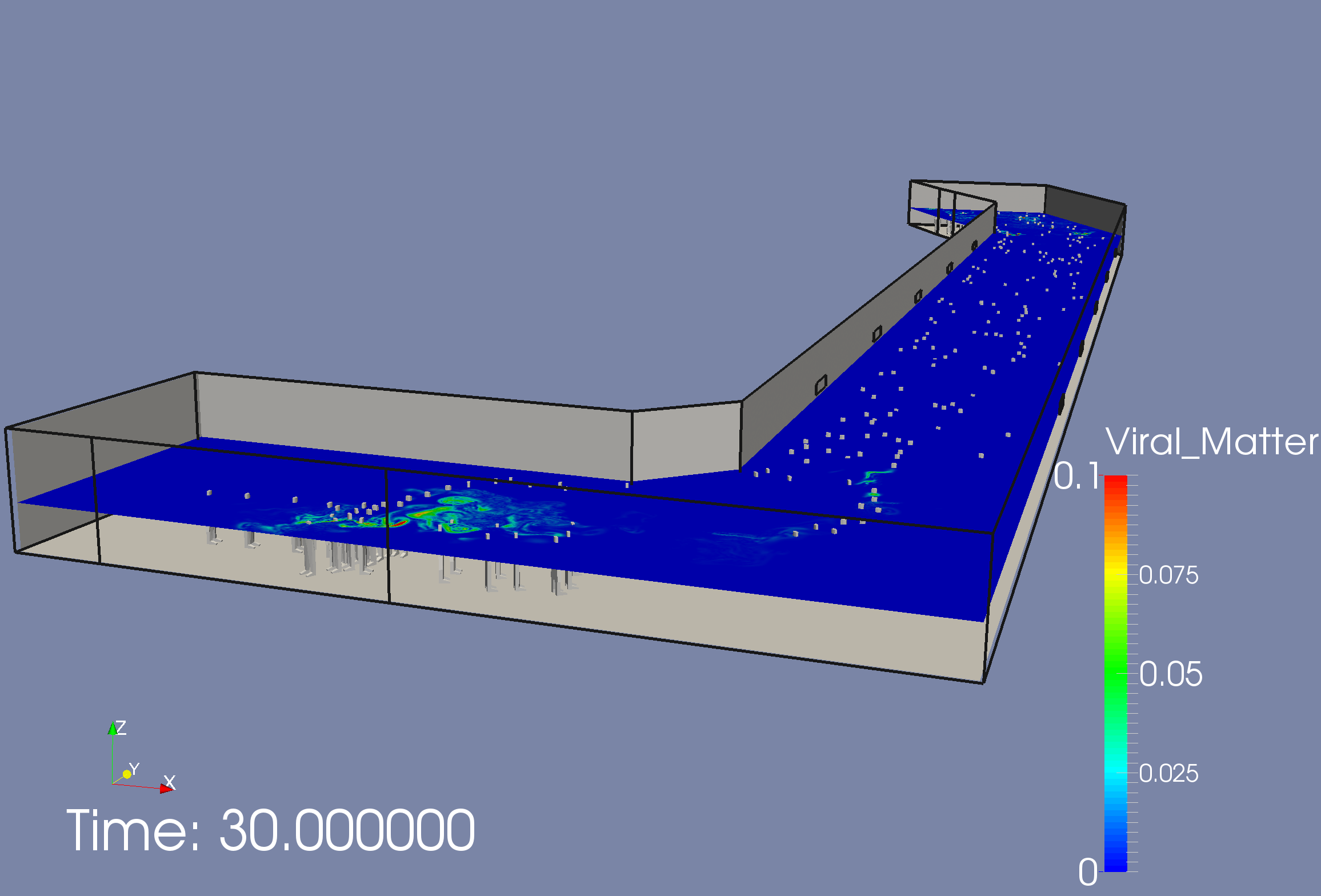}
	\caption{Long Passage: Velocity and Viral Load 
at $t=30.00~sec$}
	\label{f:pass_2}
\end{figure}

\begin{figure}
	\includegraphics[width=6.0cm]{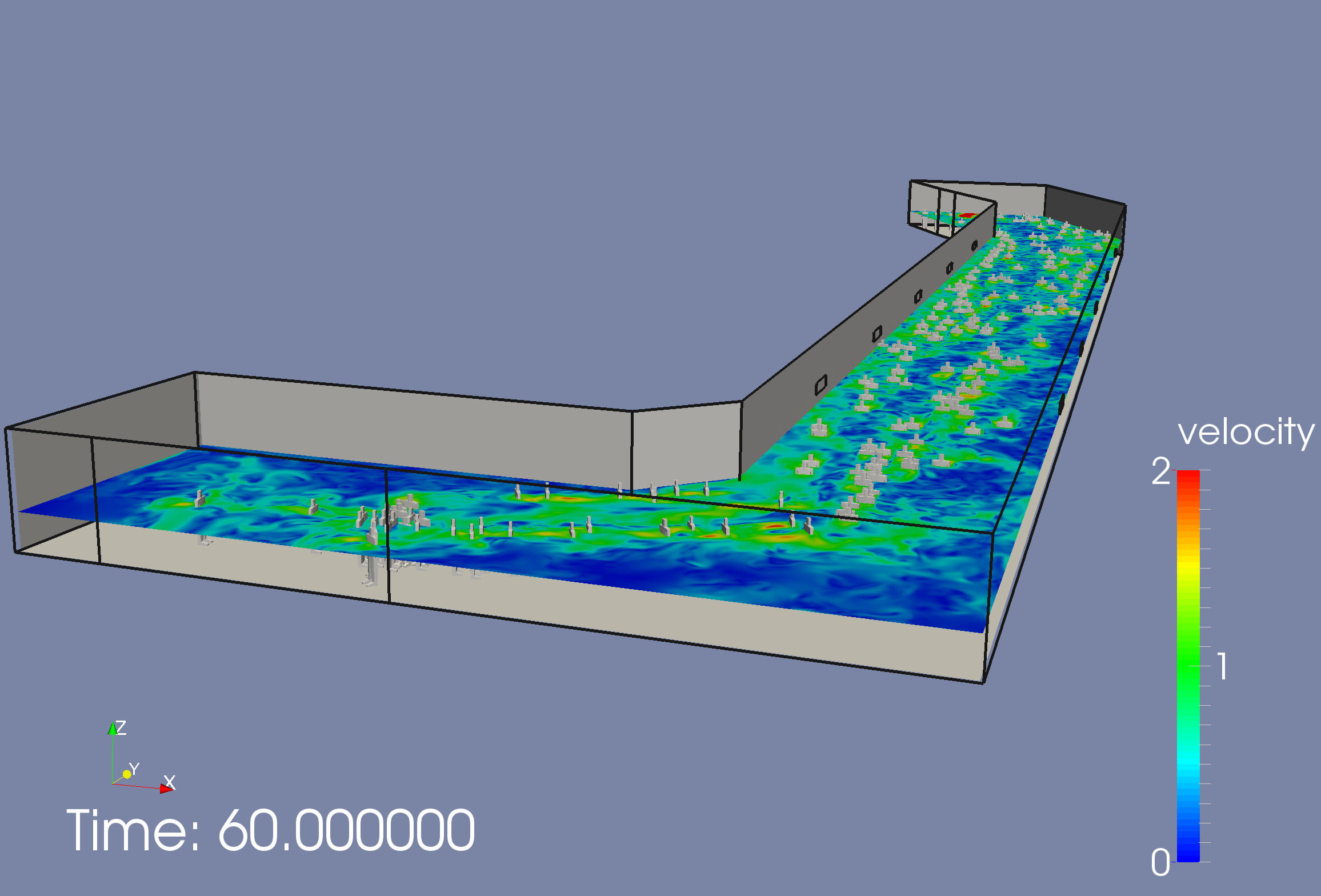}
	\includegraphics[width=6.0cm]{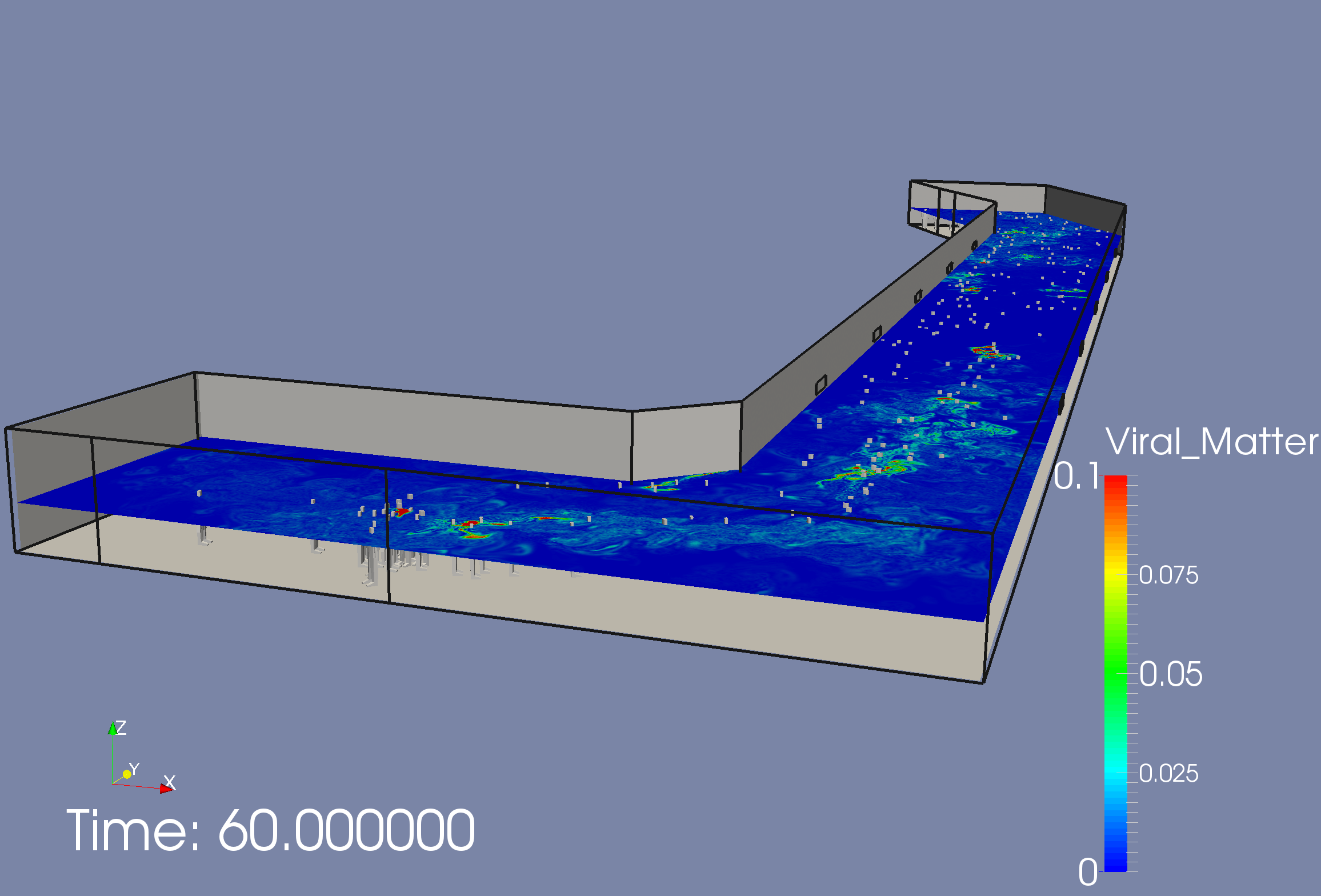}
	\caption{Long Passage: Velocity and Viral Load 
at $t=60.00~sec$}
	\label{f:pass_3}
\end{figure}

\begin{figure}
	\includegraphics[width=6.0cm]{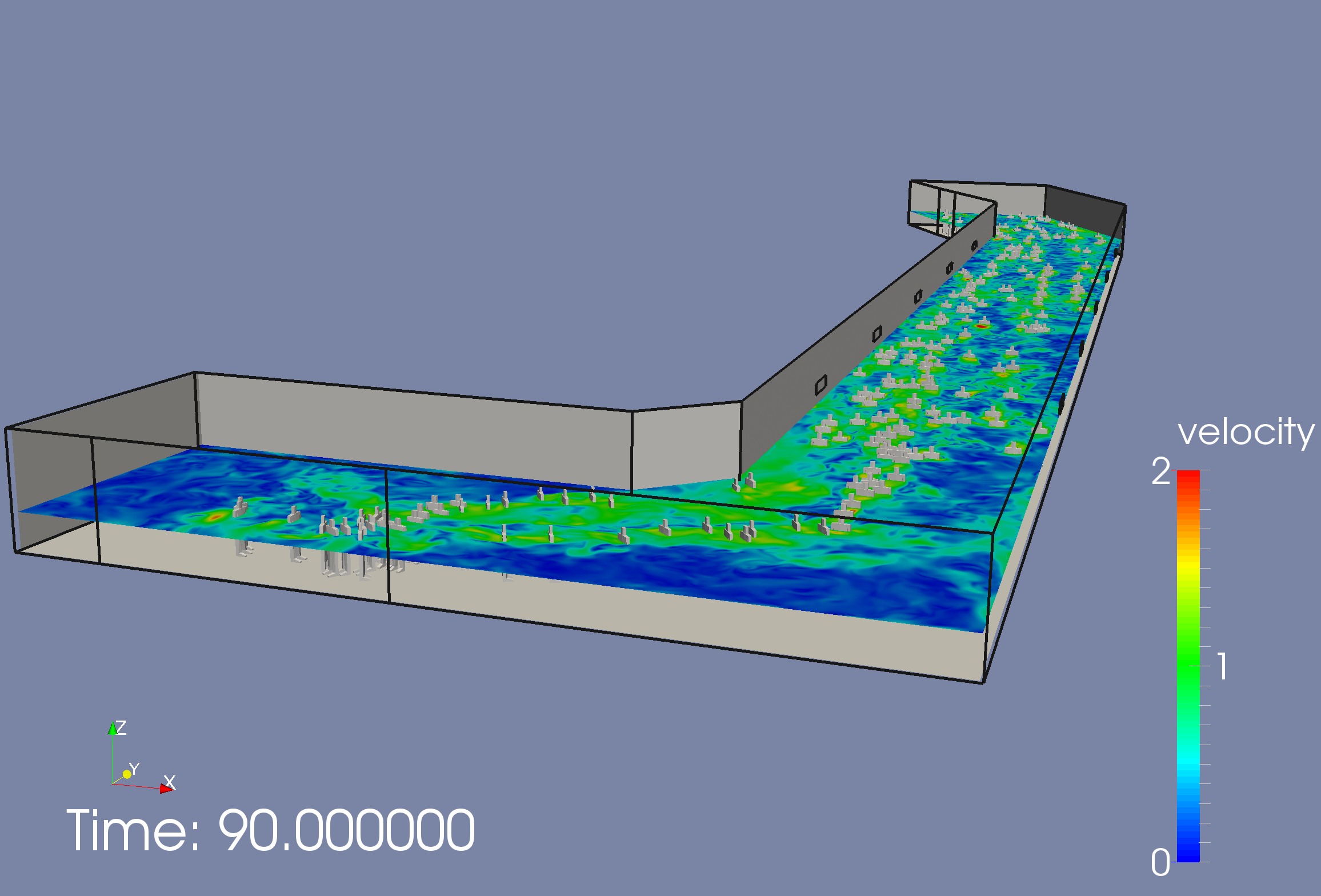}
	\includegraphics[width=6.0cm]{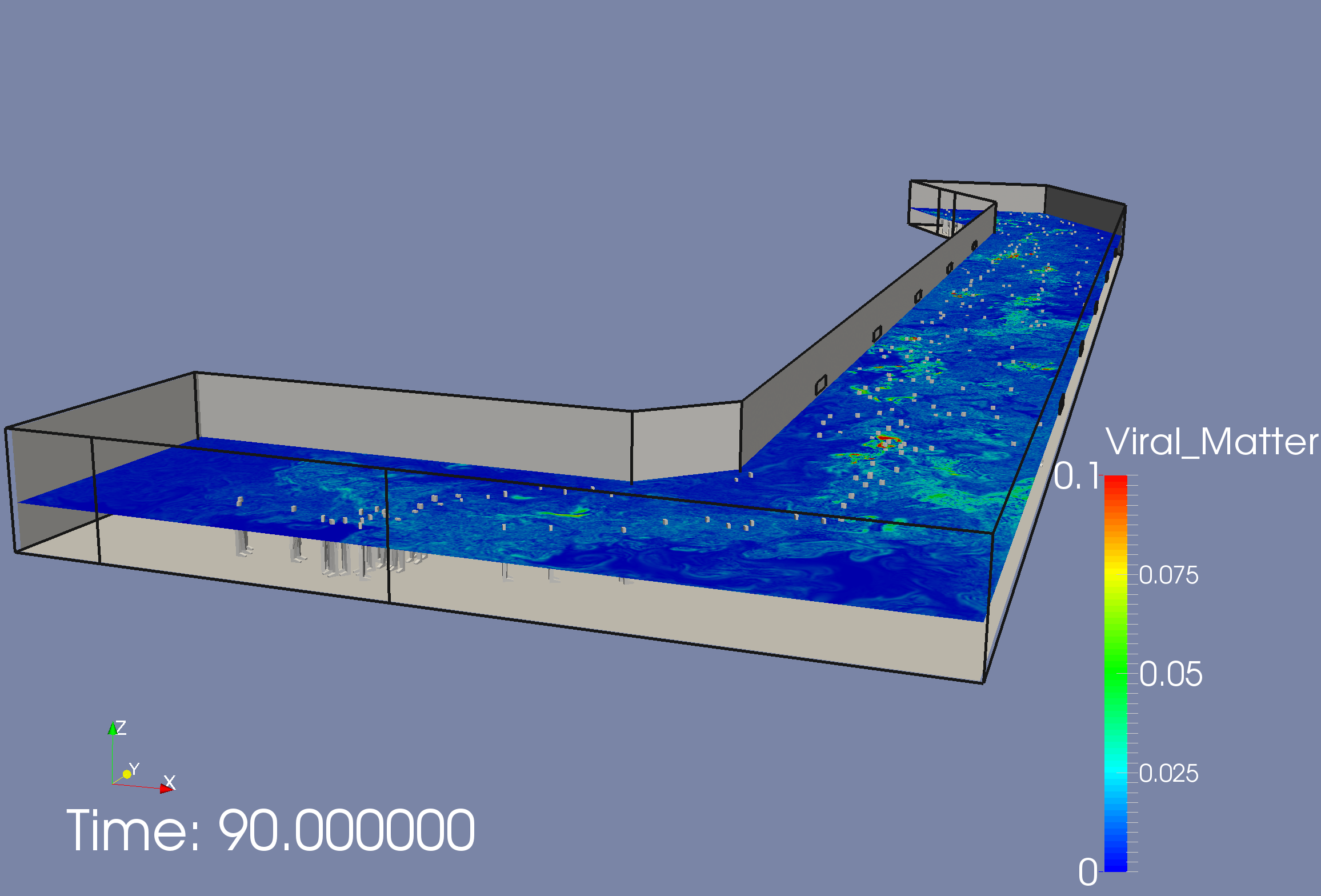}
	\caption{Long Passage: Velocity and Viral Load 
at $t=90.00~sec$}
	\label{f:pass_4}
\end{figure}

\begin{figure}
	\includegraphics[width=6.0cm]{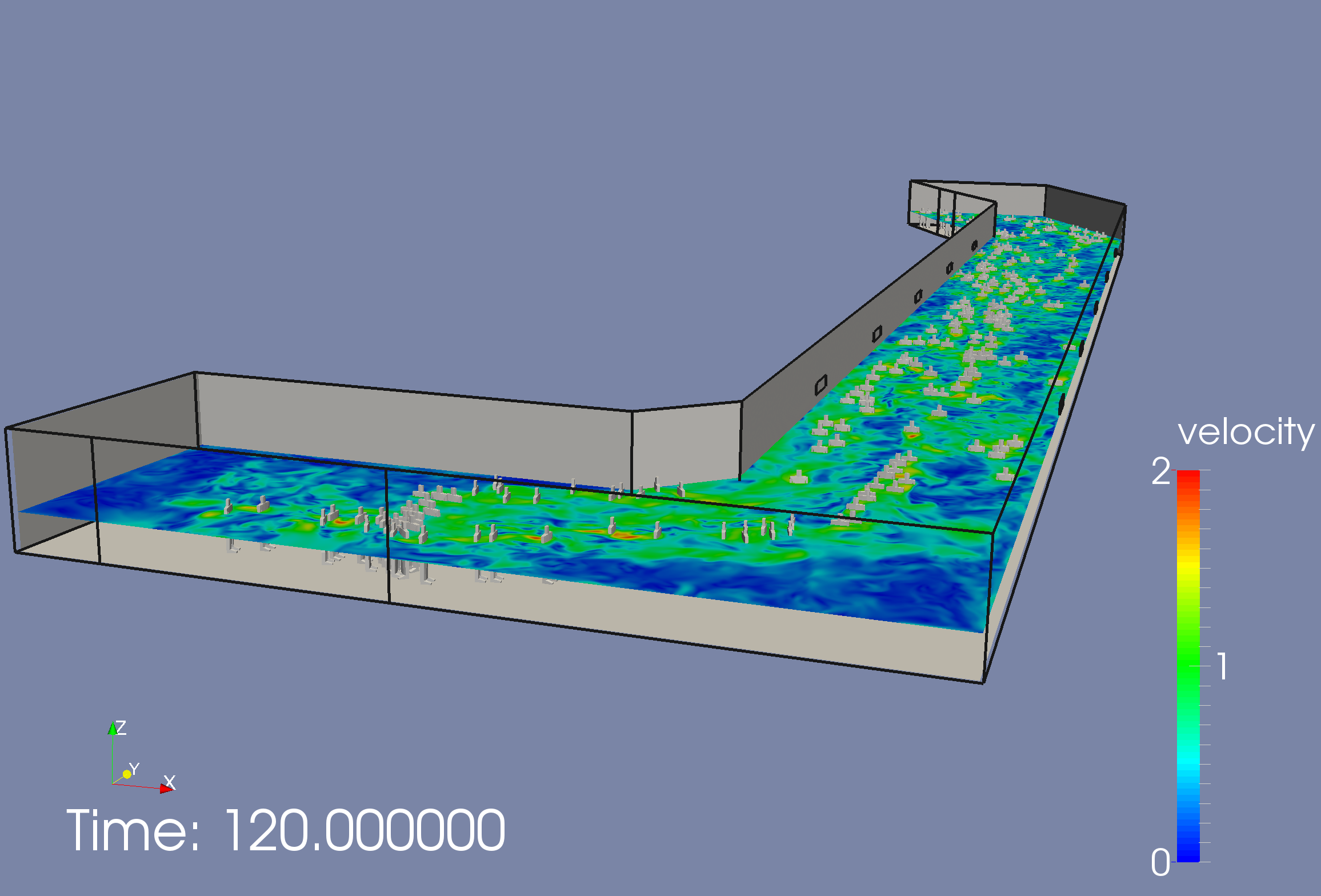}
	\includegraphics[width=6.0cm]{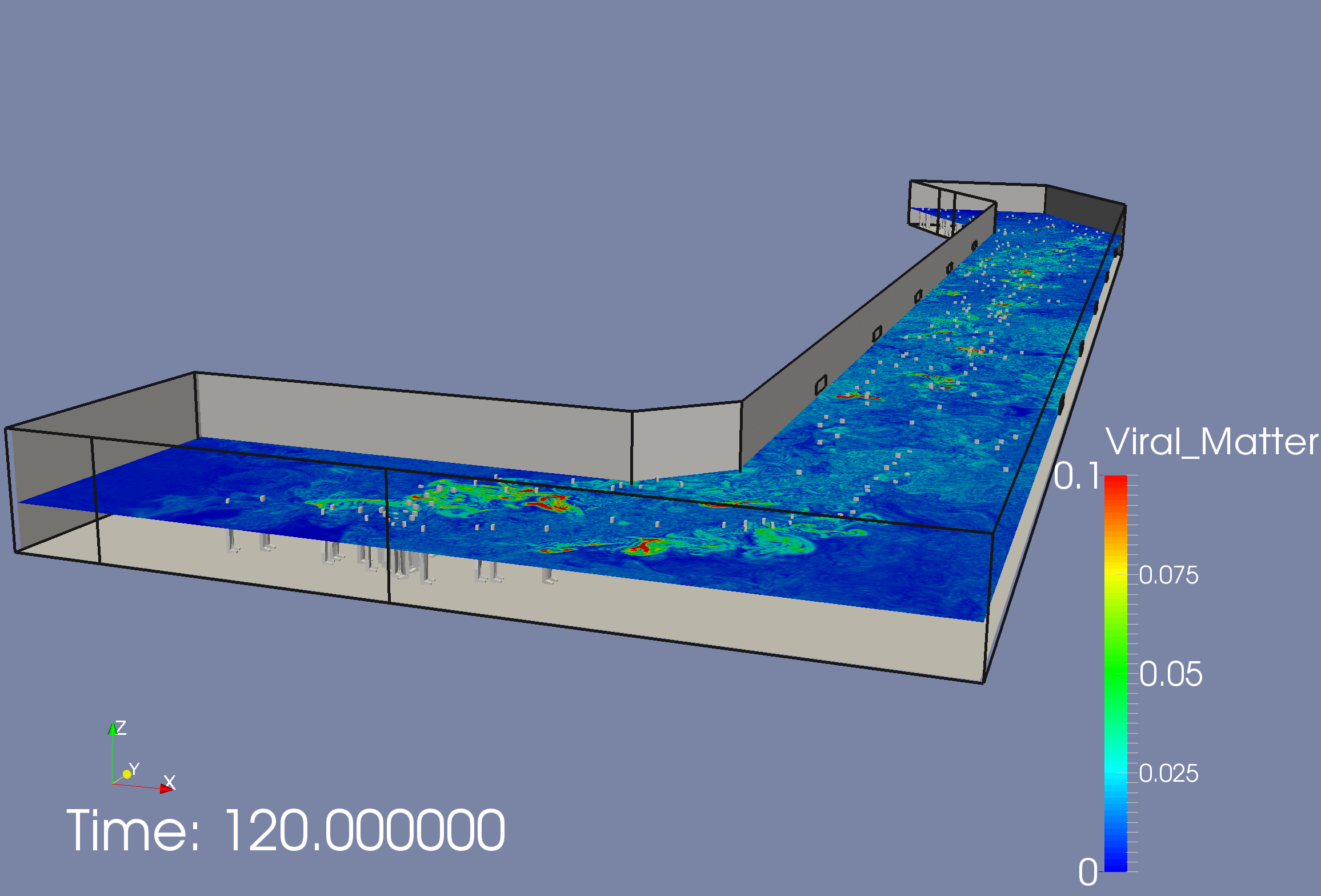}
	\caption{Long Passage: Velocity and Viral Load
at $t=120.00~sec$}
	\label{f:pass_5}
\end{figure}

\begin{figure}
	\includegraphics[width=9.0cm]{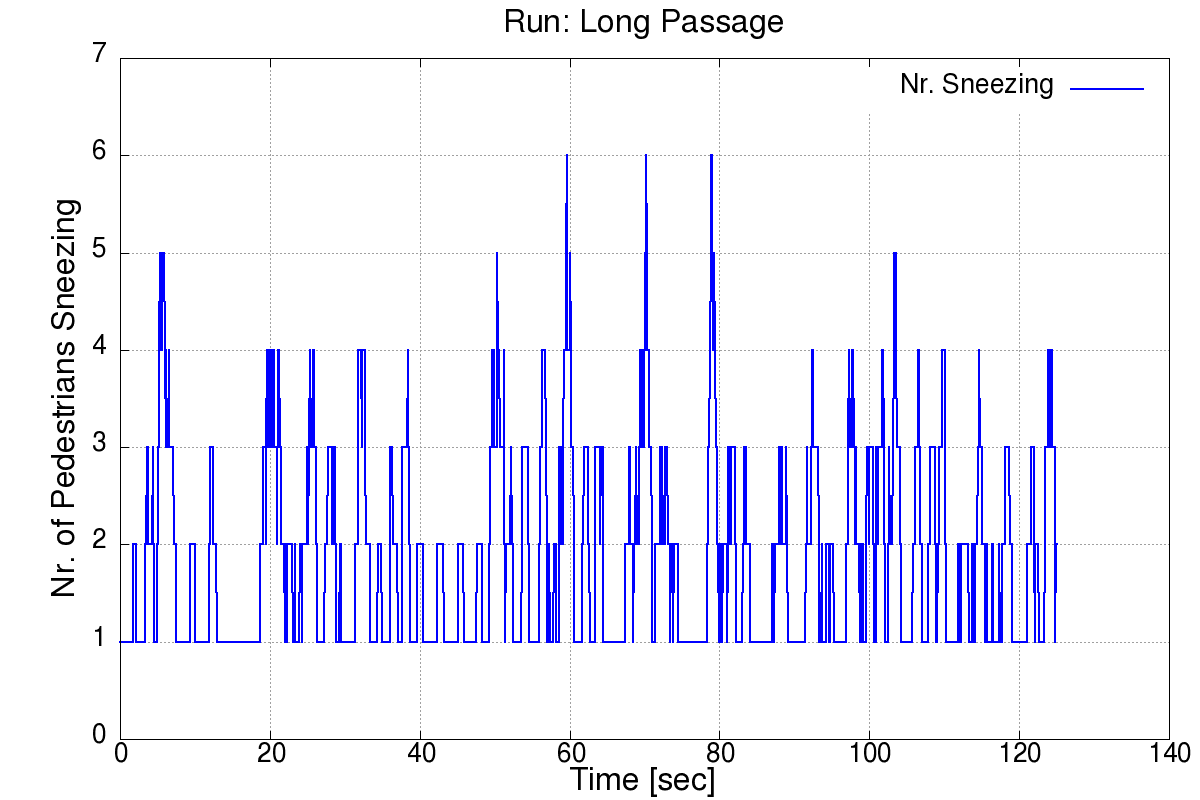}
	\caption{Long Passage: Number of Pedestrians Sneezing}
	\label{f:pass_6}	
\end{figure}

\begin{figure}
	\includegraphics[width=9.0cm]{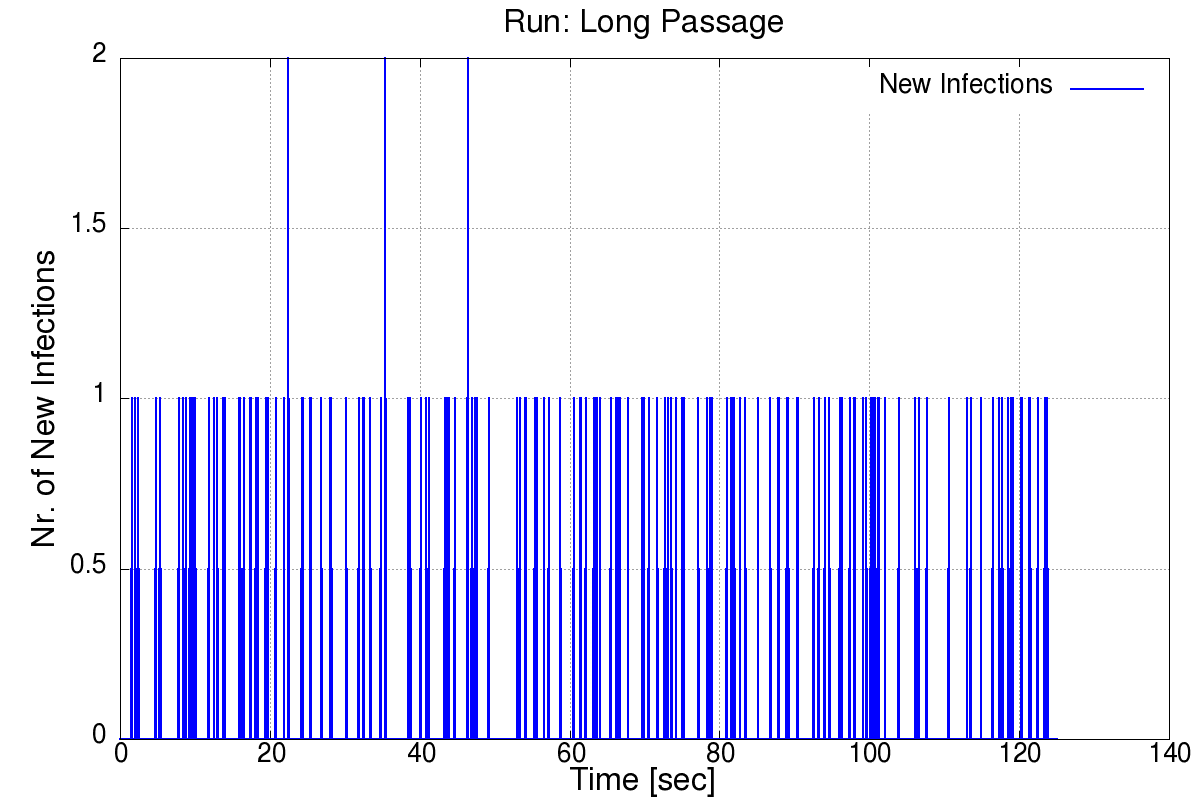}
	\caption{Long Passage: Number of New Infections}
	\label{f:pass_7}	
\end{figure}

\begin{figure}
	\includegraphics[width=9.0cm]{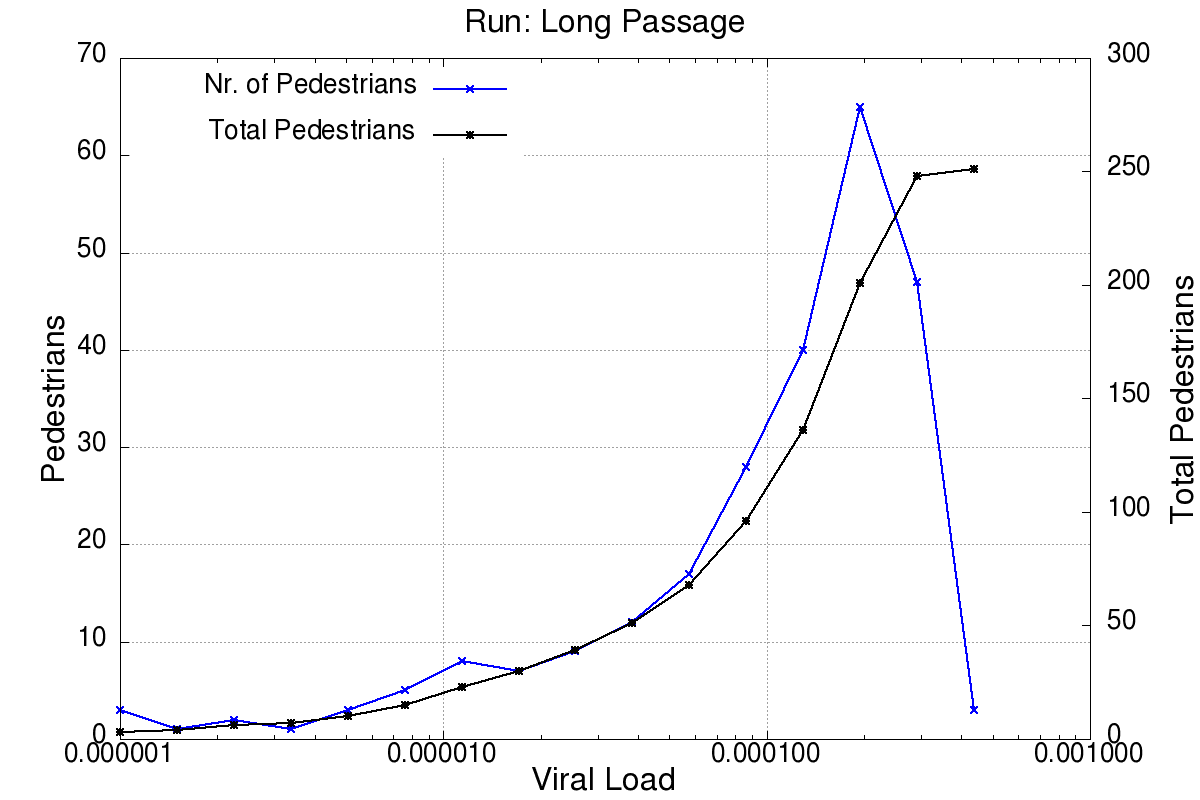}
	\caption{Long Passage: Viral Load Statistics at $t=120.00~sec$. 
		     Blue Curve: Number of Pedestrians That Inhaled a Given Viral Load. 
		     Black Curve: Cumulative Number of Pedestrians}
	\label{f:pass_8}		     
\end{figure}

This case considers a typical rush hour scenario in a train station.
The geometry and initial conditions are shown in Figure~\ref{f:pass_1}. The
central part of the passage is 82~m long, 10~m wide and 4~m high.
Several open windows are located along the passage. The inflow
velocity of the windows on the right side was set to $v_{in}=1~m/sec$.
A steady stream of pedestrians enters through both extremes and
then moves towards the exit on the other side. The pedestrian flux
is $2~p/sec$, i.e. $1~p/sec$ per entrance. This leads to an average
number of approximately 250~pedestrians in the simulation domain 
at any given time during the simulaiton. At the beginning, all
pedestrians are considered healthy. Thereafter, 70\% of the pedestrians
entering are assigned as healthy, 10\% as infected but
asymptomatic (i.e. not infecting) and 20\% as infected and infecting.
For the infecting pedestrians, the average time between sneezing
is set to $\Delta t_s=30~sec$ with a variation of 
$\delta t_s=5~sec$ and a sneeze duration of $T_s=1~sec$.
The mesh size in the region of the pedestrians was set to
$h=0.0625~m$, which implies approximately 8~flow points per pedestrian
cross-section of $d_p=0.50~m$. Not very accurate, but sufficient in
order to obtain the vortical motion due to pedestrian movement.
The total mesh size was in excess of $1.2 \cdot 10^8$~elements.
The two minutes of simulated time took approximately 8~hours using
1,024~cores. The results obtained are shown in Figures~\ref{f:pass_2}--\ref{f:pass_5}.
Note the gradual buildup of viral matter, even though fresh air is
streaming in through the windows. The simulation clearly shows the
need for better ventilation: after a while, almost all pedestrians
are wading through a `soup of viral matter'. The number of pedestrians
sneezing at any given moment has been recorded in Figure~\ref{f:pass_6}.
For a given (ableit arbitrary) threshold of $10^{-4}$ units of
viral matter, the number of new infections is shown in Figure~\ref{f:pass_7}.
The amount of viral load inhaled by the pedestrians in the 
simulation at time $T=120~sec$ may be discerned in Figure~\ref{f:pass_8}.

\newpage

\section{Conclusions and Outlook}
\label{sec:conclusions}
A deterministic pathogen transmission model based on high-fidelity physics
has been developed. The model combines computational fluid dynamics and
computational crowd dynamics in order to be able to provide accurate
tracing of viral matter that is exhaled, transmitted and inhaled via 
aerosols. The examples
shown indicate that even with modest computing resources, the propagation
and transmission of viral matter can be simulated for relatively large 
areas with thousands of square meters, hundreds of pedestrians and several
minutes of physical time. The results obtained and insights gained from 
these simulations can be used to inform global pandemic propagation models, 
increasing substantially their accuracy. \\
As with any technology, further advances are clearly possible.
The list is long, and we just mention:
\begin{itemize}
\item[-] Improved knowledge of the amount of virons in the droplets
exhaled by infecting individuals;
\item[-] Improved knowledge of the infectious dose required to 
trigger infection/illness;
\item[-] Improved boundary conditions for HVAC exits; and
\item[-] Improved modeling of particle retention and movement
through cloths (e.g. for masks); and
\item[-] Improved knowledge of the effectivity of filters in HVAC
systems that recirculate a percentage of the air in rooms.
\end{itemize}
\par \noi
Furthermore, even though the basic physical phenomena and the
partial and ordinary differential equations describing them
have been known for over a century, and solvers have advanced 
considerably over the last four decades, a vigorous experimental 
program is needed to complement and validate the numerical methods, 
and to establish firm `best practice' guidelines. 


\bibliographystyle{aiaa}

\begin{thebibliography}{7}

\bibitem{Abu20} M. Abuhegazy, K. Talaat, O. Anderoglu and
S.V. Poroseva - Numerical Investigation of Aerosol
Transport in a Classroom with Relevance to COVID-19;
{\sl Phys.\ Fluids } 32, 103311 (2020); 
https://doi.org/10.1063/5.0029118

\bibitem{Arm07} T.W. Armstrong and C.N. Haas - A Quantitative Microbial 
Risk Assessment Model for Legionnaires' Disease: Animal Model Selection 
and Dose-Response Modeling; {\sl Risk Anal.\ }27(6):1581-1596 (2007).
https://doi.org/10.1111/j.1539-6924.2007.00990.x

\bibitem{Asa19} S. Asadi, A.S. Wexler, C.D. Cappa, S. Barreda, N.M. Bouvier
and W. Ristenpart - Aerosol Emission and Superemission During Human
Speech Increase with Voice Loudness;
{\sl Nature Scientific Reports }9, (1):2348 (2019).
www.nature.com/scientificreports/
https://doi.org/10.1038/s41598-019-38808-z

\bibitem{Asa20a} S. Asadi, A.S. Wexler, C.D. Cappa, S. Barreda, N.M. Bouvier
and W. Ristenpart - Effect of
Voicing and Articulation Manner on Aerosol Particle
Emission During Human Speech {\sl PLoS ONE} 15(1):e0227699 (2020).
https://doi.org/10.1371/journal.pone.0227699

\bibitem{Asa20b} S. Asadi, N.M. Bouvier, A.S. Wexler and W. Ristenpart -
The Coronavirus Pandemic and Aerosols: Does
COVID-19 Transmit via Expiratory Particles ?
{\sl Aerosol Science and Technology } (2020).
https://doi.org/10.1080/02786826.2020.1749229


\bibitem{Blu98} V.J. Blue and J.L. Adler - Emergent Fundamental
Pedestrian Flows from Cellular Automata Microsimulation;
{\sl Transportation Research Record }1644, 29-36 (1998).

\bibitem{Blu02} V.J. Blue and J.L. Adler - Flow Capacities from Cellular
Automata Modeling of Proportional Splits of Pedestrians by Direction;
pp.~115-22 in {\sl Pedestrian and
Evacuation Dynamics} (M. Schreckenberg and S.D. Sharma eds.),
Springer (2002).

\bibitem{Bor92} J.P. Boris, F.F. Grinstein, E.S. Oran, and R.J. Kolbe -
New Insights Into Large Eddy Simulation;
{\sl Fluid Dynamics Research }10, 199-228 (1992).

\bibitem{Bro17}
A.F. Brouwer, M.H. Weir, M.C. Eisenberg, R. Meza and J.N.S. Eisenberg -
Dose-Response Relationships for Environmentally Mediated Infectious 
Disease Transmission Models; {\sl PLoS Comput.\ Biol.\ }13(4): e1005481 
(2017). 
https://doi.org/ 10.1371/journal.pcbi.1005481.

\bibitem{Bun06} H.-J. Bungartz and
M. Sch\"afer (eds.) {\sl Fluid-Structure Interaction},
Springer Lecture Notes in
Computational Science and Engineering, Springer (2006).

\bibitem{Bus20} G. Busco, S.R. Yang, J. Seo and Y.A. Hassan -
Sneezing and Asymptomatic Virus Transmission;
{\sl Phys.\ Fluids }32, 073309 (2020)
https://doi.org/10.1063/5.0019090

\bibitem{Ceb05} J.R. Cebral and R. L\"ohner -
On the Loose Coupling of Implicit Time-Marching Codes;
{\sl AIAA}-05-1093 (2005).

\bibitem{Cam04} F. Camelli and R. L\"ohner - Assessing Maximum 
Possible Damage
for Contaminant Release Events; {\sl Engineering Computations } 21, 7,
748-760 (2004).

\bibitem{Cam04a} F. Camelli, R. L\"ohner, W.C. Sandberg and R. Ramamurti -
VLES Study of Ship Stack Gas Dynamics; {\sl AIAA}-04-0072 (2004).

\bibitem{Cam06} F. Camelli and R. L\"ohner - VLES Study of Flow
and Dispersion Patterns in Heterogeneous Urban Areas;
{\sl AIAA}-06-1419 (2006).





\bibitem{Cha09} C. Chao, M.P. Wan, L. Morawska, G. Johnson, R. Graham, 
Z. Ristovski M. Hargreaves, K. Mengersen, L. Kerrie C. Steve, 
Y. Li, X. Xie and S. Katoshevski - Characterization of
Expiration Air Jets and Droplet Size Distributions Immediately at the Mouth
Opening; {\sl J.\ of Aerosol Science }40, 2, 122-133 (2009).

\bibitem{Cou05} N. Courty and S. Musse - Simulation of Large Crowds
Including Gaseous Phenomena; pp.206-212 in {\sl Proc.\ IEEE
Computer Graphics International 2005}, New York, June (2005).

\bibitem{Cur12} S. Curtis and D. Manocha - Pedestrian Simulation Using
Geometric Reasoning in Velocity Space; 
{\sl Pedestrian and Evacuation Dynamics} (U. Weidmann, U. Kirsch
and M. Schreckenberg eds.), Springer, Heidelberg (2012).

\bibitem{Dbo20a} T. Dbouk and D. Drikakis - On Coughing and Airborne 
Droplet Transmission to Humans; {\sl Phys.\ Fluids 32}, 053310 (2020).

\bibitem{Dbo20b} T. Dbouk and D. Drikakis -
On Respiratory Droplets and Face Masks;
{\sl Phys.\ Fluids }32, 063303 (2020).
https://doi.org/10.1063/5.0015044

\bibitem{Dij02} J. Dijkstra, J. Jesurun and H. Timmermans - A Multi-Agent
Cellular Automata Model of Pedestrian Movement; pp.~173-180 in
{\sl Pedestrian and Evacuation Dynamics} (M. Schreckenberg and S.D.
Sharma eds.), Springer (2002).



\bibitem{Fra97} D.R. Franz, P.B. Jahrling, A.M. Friedlander, et al. -
Clinical Recognition and Management of Patients Exposed to
Biological Warfare Agents; {\sl JAMA }278:399e411 (1997).

\bibitem{Fru71} J.J. Fruin - {\sl Pedestrian Planning and Design};
Metropolitan Association of Urban Designers and Environmental
Planners, New York (1971).

\bibitem{Fur99} C. Fureby and F. Grinstein - Monotonically Integrated
Large Eddy Simulation of Free Shear Flows; {\sl AIAA J.\ }37, 5,
544-556 (1999).

\bibitem{Gim21} J.M. Gimenez, S.R. Idelsohn, E. O\~nate and
R. L\"ohner - A Multiscale Approach for the Numerical Simulation of 
Turbulent Flows with Droplets; {\sl Archives of Computational 
Methods in Engineering }28, 4185-4204 (2021.
https://doi.org/10.1007/s11831-021-09614-6

\bibitem{Gri02} F.F. Grinstein and C. Fureby - Recent Progress on MILES
for High-Reynolds-Number Flows; {\sl J.\ Fluids Engineering }124,
848-861 (2002).


\bibitem{Gup09} J.K. Gupta, C-H. Lin and Q. Chen -
Flow Dynamics and Characterization of a Cough; 
{\sl Indoor Air }19, 517-525 (2009).

\bibitem{Gup10} J.K. Gupta, C-H. Lin and Q. Chen -
Characterizing Exhaled Airflow from Breathing and Talking; 
{\sl Indoor Air } 20, 31-39 (2010).


\bibitem{Gup12} J.K. Gupta, C-H. Lin and Q. Chen -
Risk Assessment of Airborne Infectious Diseases in Aircraft Cabins;
{\sl Indoor Air }22(5), 388-395 (2012).

\bibitem{Gup11b} J.K. Gupta, C-H. Lin and Q. Chen -
Transport of Expiratory Droplets in an Aircraft Cabin;
{\sl Indoor Air } 21, 3-11 (2011).

\bibitem{Guy09}  S.J. Guy, J. Chhugani, C. Kim, N. Satish, M. Lin,
D. Manocha and P. Dubey - ClearPath: Highly Parallel
Collision Avoidance for Multi-Agent Simulation; pp. 177-187 in {\sl
Proc.\ ACM SIGGRAPH/ Eurographics Symposium
on Computer Animation} (D. Fellner and S. Spencer eds),
Association of Computing Machinery, New York (2009).

\bibitem{Guy10} S.J. Guy, J. Chhugani, S. Curtis, P. Dubey, M. Lin
and D. Manocha - PLEdestrians: A Least-Effort Approach
to Crowd Simulation; {\sl Eurographics/
ACM SIGGRAPH Symposium on Computer Animation},
Madrid, Spain (2010).

\bibitem{Hal12} S.K. Halloran, A.S. Wexler and W.D. Ristenpart -
A Comprehensive Breath Plume Model for Disease Transmission via 
Expiratory Aerosols; {\sl PLoS ONE }7(5):e37088 (2012). 
https://doi.org/10.1371/journal.pone.0037088

\bibitem{Hel95} D. Helbing and P. Molnar - Social Force Model for
Pedestrian Dynamics; {\sl Phys.\ Rev.\ E}, 51:4282-4286 (1995).

\bibitem{Hel97} D. Helbing and P. Molnar - Self-Organization Phenomena
in Pedestrian Crowds;  569-577 in {\sl Self-Organization of Complex
Structures: From Individual to Collective Dynamics}
(F. Schweitzer (Ed.), London: Gordon and Breach (1997).

\bibitem{Hel02} D. Helbing, I.J. Farkas, P. Moln\'ar and T. Vicsek -
Simulation of Pedestrian Crowds in Normal and Evacuation Situations;
pp.~21-58 in {\sl Pedestrian and Evacuation Dynamics}
(M. Schreckenberg and S.D. Sharma eds.), Springer (2002).



\bibitem{Hug03} R.L. Hughes - The Flow of Human Crowds;
{\sl Annual Review of Fluid Mechanics }35, 169-182 (2003).

\bibitem{Ide19} S.R. Idelsohn, N. Nigro, A. Larreteguy, J.M. Gimenez
and P. Ryshakov - A Pseudo-DNS Method for the Simulation of
Incompressible Fluid Flows with Instabilities at Different Scales;
{\sl  Int. J. Comp.\ Particle Mechanics} (2019).
https://doi.org/10.1007/s40571-019-00264-x

\bibitem{Ide21} S.R. Idelsohn, J.M. Gimenez, N.M. Nigro and
E. O\~nate - The Pseudo-Direct Numerical Simulation Method for 
Multi-Scale Problems in Mechanics; {\sl Comput.\ Meth.\ Appl.\ Mech.\
Engrg.\ }380, 113774 (2021).
https://doi.org/10.1016/j.cma.2021.113774

\bibitem{Ip07} M. Ip, J.W. Tang, D.S.C. Hui, A.L.N. Wong, M.T.V. Chan, 
G.M. Joynt, A.T.P. So, S.D. Hall, P.K.S. Chan and J.J.Y. Sung -
Airflow and Droplet Spreading Around Oxygen Masks: A Simulation Model
for Infection Control Research;
{\sl AJIC } 35, 10, 684-689 (2007).

\bibitem{Ise14a}
M. Isenhour and R. L\"ohner - Verification of a Pedestrian
Simulation Tool Using the NIST Recommended Test Cases;
{\sl The Conference in Pedestrian and Evacuation Dynamics
2014 (PED2014), Transportation Research Procedia }2, 237-245 (2014).

\bibitem{Ise14b}
M. Isenhour and R. L\"ohner - Verification of a Pedestrian
Simulation Tool Using the NIST Stairwell Evacuation Data;
{\sl The Conference in Pedestrian and Evacuation Dynamics
2014 (PED2014), Transportation Research Procedia }2, 739-744 (2014).

\bibitem{Ise16} M. Isenhour - Simulating Occupant
Response to Emergency Situations; {\sl PhD Thesis},
George Mason University, Fairfax, VA (2016).

\bibitem{Ise16a} M. Isenhour and R. L\"ohner -
Validation Data from the Evacuation of a Student Center;
pp. 472-479 in
{\sl Proc. Pedestrian and Evacuation Dynamics 2016 (PED~2016)},
(W. Song, J. Ma and L. Fu eds.), University of Science and
Technology Press, Hefei, China, Oct 17-21 (2016).

\bibitem{Ise16b} M. Isenhour and R. L\"ohner -
Pedestrian Speed on Stairs: A Mathematical Model Based on
Empirical Analysis for Use in Computer Simulations;
pp. 529-533 in
{\sl Proc. Pedestrian and Evacuation Dynamics 2016 (PED~2016)},
(W. Song, J. Ma and L. Fu eds.), University of Science and
Technology Press, Hefei, China, Oct 17-21 (2016).

\bibitem{Jam81}
A. Jameson, W. Schmidt and E. Turkel - Numerical Solution
of the Euler Equations by Finite Volume Methods using Runge-Kutta
Time-Stepping Schemes; {\sl AIAA}-81-1259 (1981).

\bibitem{Joh11} G.R. Johnson, L. Morawska, Z.D. Ristovski,
M. Hargreaves, K. Mengersen,
C.Y.H. Chao, M.P. Wan, Y. Li , X. Xie, D. Katoshevski, S. Corbette -
Modality of Human Expired Aerosol Size Distributions
{\sl J.\ of Aerosol Science }42, 839-851 (2011).


\bibitem{Kar10} T. Karmakharm, P. Richmond and D.M. Romano -
Agent-based Large Scale Simulation of Pedestrians With
Adaptive Realistic Navigation Vector Fields;
{\sl EG UK Theory and Practice of Computer Graphics 2010}
(J. Collomosse and I. Grimstead eds.) (2010).

\bibitem{Kes02} A. Kessel, H. Kl\"upfel, J. Wahle and M. Schreckenberg -
Microscopic Simulation of Pedestrian Crowd Motion; pp.~193-202 in
{\sl Pedestrian and Evacuation Dynamics} (M. Schreckenberg and S.D.
Sharma eds.), Springer (2002).



\bibitem{Klu03} H.L. Kl\"upfel - A Cellular Automation Model for Crowd
Movement and Egress Simulation; {\sl Ph.D. Dissertation: Falkut\"at }4,
Univ. Duisburg-Essen (2003).

\bibitem{Lak05} T.I. Lakoba, D.J. Kaup and N.M. Finkelstein - Modifications
of the Helbing-Moln\'ar-Farkas-Vicsek Social Force Model for Pedestrian
Evolution; {\sl Simulation }81, 339 (2005).

\bibitem{Lan06} P.A. Langston, R. Masling and B.N. Asmar -
Crowd Dynamics Discrete Element Multi-Circle Model;
{\sl Safety Science }44, 395-417 (2006).


\bibitem{Li20} H. Li, F.Y. Leong, G. Xu, Z. Ge, Ch.W. Kang and K.H. Lim -
Dispersion of Evaporating Cough Droplets in Tropical Outdoor 
Environment; {\sl Phys.\ Fluids }32, 113301 (2020)
https://doi.org/10.1063/5.0026360

\bibitem{Lin10}
W.G. Lindsley, F.M. Blachere, R.E. Thewlis, A. Vishnu, K.A. Davis, 
G. Cao, et al.\ - Measurements of Airborne Influenza Virus in Aerosol 
Particles from Human Coughs; {\sl PLoS ONE }5:e15100 (2010).

\bibitem{Lin12}
W.G. Lindsley, T.A. Pearce, J.B. Hudnall, K.A. Davis, S.M. Davis, 
M.A. Fisher, et al.\ - Quantity and Size Distribution of Cough-Generated 
Aerosol Particles Produced by Influenza Patients During and After 
Illness; {\sl J.\ Occup.\ Environ.\ Hyg.\ } 9, 443-9. (2012).

\bibitem{Loh95} R. L\"ohner, C. Yang, J. Cebral,
J.D. Baum, H. Luo, D. Pelessone and C. Charman - Fluid-Structure
Interaction Using a Loose Coupling Algorithm
and Adaptive Unstructured Grids; {\sl AIAA}-95-2259 [Invited] (1995). -

\bibitem{Loh01} R. L\"ohner, Chi Yang, J. Cebral, O. Soto, F. Camelli,
J.D. Baum, H. Luo, E. Mestreau, D. Sharov, R. Ramamurti, W. Sandberg
and Ch. Oh - Advances in FEFLO; {\sl AIAA}-01-0592 (2001).

\bibitem{Loh02} R. L\"ohner, Chi Yang, J. Cebral, O. Soto, F. Camelli,
J.D. Baum, H. Luo, E. Mestreau and D. Sharov - Advances in FEFLO;
{\sl AIAA}-02-1024 (2002).

\bibitem{Loh04} R. L\"ohner - Multistage Explicit Advective Prediction for
Projection-Type Incompressible Flow Solvers;
{\sl J.\ Comp.\ Phys.\ }195, 143-152 (2004).

\bibitem{Loh06} R. L\"ohner, Chi Yang,
J.R. Cebral,
F. Camelli, O. Soto and J. Waltz - Improving the Speed and Accuracy
of Projection-Type Incompressible Flow Solvers;
{\sl Comp.\ Meth.\ Appl.\ Mech.\ Eng.\ }195, 23-24, 3087-3109
(2006).

\bibitem{Loh08} R. L\"ohner - {\sl Applied CFD Techniques, 
Second Edition}; J. Wiley \& Sons (2008).

\bibitem{Loh08b} R. L\"ohner, J.R. Cebral, F.F. Camelli,
S. Appanaboyina, J.D. Baum, E.L. Mestreau and O. Soto - Adaptive
Embedded and Immersed Unstructured Grid Techniques;
{\sl Comp.\ Meth.\ Appl.\ Mech.\ Eng.\ }197, 2173-2197 (2008).

\bibitem{Loh10} R. L\"ohner - On the Modeling of Pedestrian Motion;
{\sl Appl. Math. Modelling }34, 2, 366-382 (2010).

\bibitem{Loh11}
R. L\"ohner - Coupling Several CFD and CSD Codes
in One Application; pp.\ 1 - 16 in Special Edition
{\sl Int.\ J.\ of Multiphysics } (2011).

\bibitem{Loh13}  R. L\"ohner, F. Camelli, J.D. Baum, O.A. Soto
and F. Togashi - Advances in FEFLO; {\sl AIAA}-2013-0373 (2013).

\bibitem{Loh14} R. L\"ohner, F. Camelli, J.D. Baum, F. Togashi and
O. Soto - On Mesh-Particle Techniques; {\sl Comp.\ Part.\ Mech.\ }1,
199-209 (2014).

\bibitem{Loh16} R. L\"ohner, M. Baqui, E. Haug and B. Muhamad - Real-Time
Micro-Modelling of a Million Pedestrians; {\sl Engineering Computations }
33, 1, 217-237 (2016).

\bibitem{Loh16a} R. L\"ohner and F. Camelli -
Tightly Coupled Computational Fluid and Crowd Dynamics; pp.
505-509 in
{\sl Proc. Pedestrian and Evacuation Dynamics 2016 (PED~2016)},
(W. Song, J. Ma and L. Fu eds.), University of Science and
Technology Press, Hefei, China, Oct 17-21 (2016).

\bibitem{Loh20} R. L\"ohner, H. Antil, S. Idelsohn and E. O\~nate -
Detailed Simulation of Viral Propagation in the Built
Environment; {\sl Computational Mechanics }66, 1093-1107 (2020).
https://doi.org/10.1007/s00466-020-01881-7 

\bibitem{Loh21} R. L\"ohner and H. Antil -
High Fidelity Modeling of Aerosol Pathogen Propagation
in Built Environments With Moving Pedestrians;
{\sl Int.\ J.\ Num.\ Meth.\ Biomed.\ Engng.\ } 37, 3 (2021);e3428.
https://doi.org/10.1002/cnm.3428

\bibitem{Loh21a} R. L\"ohner, H. Antil, A. Srinivasan, S. Idelsohn
and E. O\~nate -
High-Fidelity Simulation of Pathogen Propagation, Transmission
and Mitigation in the Built Environment;
{\sl Archives of Computational Methods in Engineering } 28, 6, 
4237-4262 (2021).
https://doi.org/10.1007/s11831-021-09606-6

\bibitem{Lou67} R.G. Loudon and R.M. Roberts - Droplet Expulsion 
from the Respiratory Tract; {\sl Am.\ Rev.\ Respir.\ Dis.\ } 95, 
3, 435-442 (1967).

\bibitem{Mil13}
D.K. Milton, M.P. Fabian, B.J. Cowling, M.L. Grantham, J.J. McDevitt -
Influenza Virus Aerosols in Human Exhaled Breath: Particle Size, 
Culturability, and Effect of Surgical Masks; {\sl PLoS Pathog.\ } 9:e1003205
(2013).

\bibitem{Nam17} S. Namilae, P. Derjany, A. Mubayi, M. Scotch and 
A. Srinivasan - Multiscale Model for Pedestrian and Infection Dynamics 
During Air Travel; {\sl Physical Review E }95(5), 052320 (2017).

\bibitem{Pel06} N. Pelechano and N.I. Badler - Modeling Crowd
and Trained Leader Behavior During Building Evacuation;
{\sl IEEE Computer Graphics and Applications }26 (6): 80-86 (2006).

\bibitem{Pel08} N. Pelechano, J. Allbeck and N.I. Badler - Virtual
Crowds: Methods, Simulation and Control; Morgan \& Claypool,
San Rafael, CA (2008).

\bibitem{Pre71} W.M. Predtetschenski and A.I. Milinski -
{\sl Personenstr\"ome in Geb\"auden - Berechnungsmethoden f\"ur die
Projektierung}; Verlaggesellschaft Rudolf M\"uller,
K\"oln-Braunsfeld (1971).

\bibitem{Qui03} M.J. Quinn, R.A. Metoyer and K. Hunter-Zaworski -
Parallel Implementation of the Social Forces Model; pp. 63-74 in
{\sl Proc.\ 2nd Int.\ Conf.\ in Pedestrian and Evacuation Dynamics} (2003).

\bibitem{Ram96} R. Ramamurti and
R. L\"ohner - A Parallel Implicit
Incompressible Flow Solver Using Unstructured Meshes;
{\sl Computers and Fluids } 5, 119-132 (1996).

\bibitem{Ram99} R. Ramamurti, W.C. Sandberg and R. L\"ohner - Computation of
Unsteady Flow Past Deforming Geometries; {\sl Int.\ J.\ Comp.\ Fluid Dyn.\ },
83-99 (1999).

\bibitem{Scha02} A. Schadschneider - Cellular Automaton Approach
to Pedestrian Dynamics - Theory; pp.~75-86 in
{\sl Pedestrian and Evacuation Dynamics}
(M. Schreckenberg and S.D. Sharma eds.), Springer (2002).


\bibitem{Sch08} M. Sch\"afer
and S. Turek (eds.) - {\sl Proc. Int.\ Workshop on
Fluid-Structure Interaction:
Theory, Numerics and Applications}, Herrsching (Munich), Germany,
Sept. 29 - Oct 1 (2008).

\bibitem{Sud07} A. Sud, R. Gayle, E. Andersen, S. Guy, Ming Lin and
D. Manocha - Real-time Navigation of Independent Agents Using
Adaptive Roadmaps; {\sl ACM Symposium on Virtual Reality Software
and Technology} (2007).

\bibitem{Sze10}
G.N. Sze To and C.Y. Chao - Review and Comparison Between the 
Wells-Riley and Dose-Response Approaches to Risk Assessment of 
Infectious Respiratory Diseases; {\sl Indoor Air }20(1):2-16 (2010).
doi:10.1111/j.1600-0668.2009.00621.x.

\bibitem{Tan06} J.W. Tang, Y. Li, I. Eames, P.K.S. Chan and G.L. Ridgway
Factors Involved in the Aerosol Transmission of Infection and 
Control of Ventilation in Healthcare Premises;
{\sl J.\ of Hospital Infection }64, 100-114 (2006).

\bibitem{Tan11} J.W. Tang, C.J. Noakes, P.V. Nielsen, I. Eames, 
A. Nicolle, Y. Li and
G.S. Settles - Observing and Quantifying Airflows in the Infection 
Control of Aerosol- and Airborne-Transmitted Diseases: An Overview of 
Approaches; {\sl J.\ of Hospital Infection }77 213-222 (2011).

\bibitem{Tan12} J.W. Tang, A.D. Nicolle, J. Pantelic, G.C. Koh, L. Wang,
M. Amin, C.A. Klettner, D.K.W. Cheong, C. Sekhar and K.W. Tham -
Airflow Dynamics of Coughing in Healthy Human Volunteers by Shadowgraph
Imaging: An Aid to Aerosol Infection Control;
{\sl PLoS ONE }7, 4: e34818 (2012).
doi:10.1371/journal.pone.0034818

\bibitem{Tan13} J.W. Tang, A.D. Nicolle, C.A. Klettner, J. Pantelic,
L. Wang, A. Bin Suhaimi, A.Y.L. Tan, G.W.X. Ong, R. Su, C. Sekhar,
D.K.W. Cheong and K.W. Tham - Airflow Dynamics of Human Jets: Sneezing and
Breathing - Potential Sources of Infectious Aerosols;
{\sl PLoS ONE }8, 4: e59970 (2013).
doi:10.1371/journal.pone.0059970

\bibitem{Tek00} K. Teknomo, Y. Takeyama and H. Inamura -
Review on Microscopic Pedestrian Simulation Model; {\sl Proc.
Japan Society of Civil Engineering Conf. } Morioka, Japan,
March (2000).

\bibitem{Teu10} P.F.M. Teunis, N. Brienen, M.E.E. Kretzschmar -
High Infectivity and Pathogenicity of Influenza A Virus Via Aerosol and
Droplet Transmission; {\sl Epidemics }2, 215-222 (2010).

\bibitem{Tha07} D. Thalmann and S.R. Musse - {\sl Crowd Simulation};
Springer-Verlag, London (2007).

\bibitem{To20} K. K.-W. To et al. - 
Temporal Profiles of Viral Load in Posterior Oropharyngeal
Saliva Samples and Serum Antibody Responses During
Infection by SARS-CoV-2: An Observational Cohort Study;
{\sl Lancet Infect.\ Dis.\ } (online) (2020).
https://doi.org/10.1016/S1473-3099(20)30196-1

\bibitem{Tor12} P.M. Torrens - Moving Agent Pedestrians Through Space
and Time; {\sl Annals of the Association of
American Geographers }102, 1, 35-66 (2012).

\bibitem{Vig10} G. Vigueras, M. Lozano, J.M. Ordun and F. Grimaldo -
A Comparative Study of Partitioning Methods for Crowd Simulations;
{\sl Applied Soft Computing} 10, 225-235 (2010).

\bibitem{Vil13} J.M. Villafruela, I. Olmedo, M. Ruiz de Adana, 
C. Mendez and P.V. Nielsen - CFD Analysis of the Human Exhalation 
Flow Using Different Boundary Conditions and Ventilation Strategies;
{\sl Building and Environment }62, 191-200 (2013).

\bibitem{Vis95} P.M. Vishton and J. E. Cutting - Wayfinding,
Displacements, and Mental Maps: Velocity Fields are Not Typically
Used to Determine One's Aimpoint; {\sl J.\ of Experimental
Psychology }21 (5): 978-995 (1995).

\bibitem{WHO20} World Health Organization - Transmission of 
SARS-CoV-2: Implications for Infection Prevention Precautions;
{\sl Scientific Brief}, July 9 (2020).
https://www.who.int/news-room/commentaries/detail/transmission-of-sars-cov-2-implications-for-infection-prevention-precautions.

\bibitem{Wei16} J. Wei, Y. Li - Airborne Spread of Infectious Agents 
in the Indoor 
Environment; {\sl American J.\ of Infection Control }44, S102-S108 (2016).
http://dx.doi.org/10.1016/j.ajic.2016.06.003

\bibitem{Wel55} W.F. Wells - {\sl Airborne Contagion and Air Hygiene.
An Ecological Study of Droplet Infections}; Cambridge University 
Press (1955).

\bibitem{Xie07} X. Xie, Y. Li, A.T.Y. Chwang, P.L. Ho, W.H. Seto -
How Far Droplets Can Move in Indoor Environments - Revisiting
the wells Evaporation-Falling Curve;
{\sl Indoor Air }17, 211-225 (2007). 
doi:10.1111/j.1600-0668.2006.00469.x

\bibitem{Zha14}
J. Zhang, D. Britto, M. Chraibi, R. L\"ohner, E. Haug and B. Gawenat -
Qualitative Validation of PEDFLOW for Description of Unidirectional
Pedestrian Dynamics;
{\sl The Conference in Pedestrian and Evacuation Dynamics
2014 (PED2014), Transportation Research Procedia }2, 733-738 (2014).

\bibitem{Zha17} Y. Zhang, G. Feng, Z. Kang, Y. Bi and Y. Cai - 
Numerical Simulation of Coughed Droplets in Conference Room;
{\sl 10th International Symposium on Heating, Ventilation and Air 
Conditioning, ISHVAC2017}, October, 19-22 Jinan, China (2017),
{\sl Procedia Engineering }205, 302-308 (2017).

\bibitem{Zoh20a} T.I. Zohdi - Modeling and 
Simulation of the Infection Zone from a Cough;
{\sl Computational Mechanics} (2020).
https://doi.org/10.1007/s00466-020-01875-5


\end{thebibliography}

\end{document}